\documentclass[thmsb,11pt]{article}
\usepackage{amsfonts}
\usepackage{amssymb}
\usepackage{graphicx}
\usepackage{amsmath}
\usepackage{geometry}

\newtheorem{theorem}{Theorem}[section]

\newtheorem{corollary}[theorem]{Corollary}

\newtheorem{definition}[theorem]{Definition}
\newtheorem{example}[theorem]{Example}

\newtheorem{lemma}[theorem]{Lemma}

\newtheorem{proposition}[theorem]{Proposition}

\input{tcilatex}
\begin{document}

\title{The significance of the contributions of congruences to the theory of
connectednesses and disconnectednesses for topological spaces and graphs.}
\author{Stefan Veldsman\bigskip \bigskip \\
\bigskip \bigskip Dedicated to the memory of Izak Broere}
\maketitle

\begin{abstract}
\noindent This is a survey of some of the consequences of the recently
introduced congruences on the theory of connectednesses (radical classes)
and disconnectednesses (semisimple classes) of graphs and topological
spaces. In particular, it is shown that the connectednesses and
disconnectednesses can be obtained as Hoehnke radicals and a connectedness
has a characterization in terms of congruences resembling the classical
characterization of its algebraic counterpart using ideals for a radical
class. But this approach has also shown that there are some unexpected
differences and surprises: an ideal-hereditary Hoehnke radical of
topological spaces or graphs need not be a Kurosh-Amitsur radical and in the
category of graphs with no loops, non-trivial connectednesses and
disconnectednesses exist, but all Hoehnke radicals degenerate.
\end{abstract}

\begin{flushleft}
\end{flushleft}

\section{Introduction}

\noindent In this survey, we report on the work providing the last piece of
the puzzle to fully establish the correspondence between the radical theory
of the algebraic structures like rings, nearrings and groups on the one
hand, and the non-algebraic categories of topological spaces and graphs on
the other. This link has been made possible by the recent definition and
development of a theory of congruences for graphs and for topological
spaces. It is then also possible to define the radical classes for
topological spaces and graphs (in these categories, such classes are called
connectednesses) as Hoehnke radicals. This opened the door to explore other
similarities and differences between the algebraic and non-algebraic
theories and we will report on some of these unexpected and surprising
difference.

\bigskip

\noindent The origins of radical theory goes back to the early twentieth
century with the work of Wedderburn on finite dimensional algebras. This was
extended to ring theory with K\"{o}the's nilradical, the Jacobson radical
and subsequently many other radicals. These, together with some developments
in group theory, led to the axiomization of the radical concept by Kurosh
and independently Amitsur for rings, groups and omega-groups in the early
fifties. In the sixties, Hoehnke used congruences to define a radical for
universal algebras; now known as a Hoehnke radical. In an environment where
a congruence is completely determined by one of its congruence classes, as
for example in $\Omega $-groups, it was then shown under which conditions a
Hoehnke radical will be a Kurosh-Amitsur radical. A second main stream that
contributed to the development and existence of general radical theory has
its origins in the torsion theory of modules. Torsion theories of modules
were defined in terms of equivalence classes of injective modules, but when
Dickson generalized torsion theories to abelain categories, this approach
was abandoned. It turns out that the torsion and torsion-free classes
correspond to the radical and semisimple classes of general radical theory
respectively. It is then interesting and rather pleasing that the hereditary
torsion theories of topological spaces have a connection with injective
topological spaces (= indiscrete spaces) as will be seen below. The third
contribution to the general radical theory is to be found in the category of
topological spaces. Connectednesses and disconnectednesses of topological
spaces were defined by Preu\ss\ as classes of spaces on which certain maps
are constant. Then Arhangel'ski\u{\i} and Wiegandt showed that these classes
corresponds to the Kurosh-Amitsur radical and semisimple classes of rings
and groups by replacing the algebraic notions involved in one of their
characterizations with their categorical suitable topologial versions. This
was followed by a theory of connectednesses and disconnectednesses for
graphs and more generally for abstract relational structures. The graphs
that we refer to in the preceding lines are graphs for which loops are
allowed. When no loops are permitted, the radical theory has many different
and unusual characteristics.

\bigskip

\noindent With such similar radical theories in so many divergent branches
of mathematics, the need arose for a common language to describe them all.
Category theory proved to be a suitable tool; initially only catering for
the radical theories from an algebraic environment, for example by \v{S}ul'ge%
\u{\i}fer [27], Suli\'{n}ski [28] and Holcombe and Walker [19]. Such
categories exclude the connectednesses and disconnectednesses of topological
spaces and graphs. A unified treatment for these two cases were given by
Fried and Wiegandt [15, 16] by considering graphs and topological spaces as
abstract relational stuctures. But this approach excluded the algebraic
cases. In [5] and [29] less stringent conditions were imposed on a general
category to make it suitable to describe the radical theory of the classical
algebraic stuctures as well as those of topological spaces and graphs.
Subsequently, the most comprehensive theory covering most known radical
theories was given by M\'{a}rki, Mlitz and Wiegandt [20] in a general
categorical setting but with an universal algebraic flavor.

\bigskip

\noindent Quite recently it has been shown that the connectednesses and
disconnectednesses of both topological spaces and graphs (i.e the
Kurosh-Amitsur radical and semisimple classes) can be obtained from Hoehnke
radicals as has been done for universal algebras using congruences. In [4]
and [30] it was shown how to define congruences on graphs and topological
spaces respectively. These congruences then lead in a natural way to
appropriate versions of the algebraic isomorphism theorems and subdirect
products. As in universal algebra, one can then define a Hoehnke radical for
graphs and topological spaces. Necessary and sufficient conditions to ensure
that the Hoehnke radical becomes a Kurosh-Amitsur radical have been
determined for algebraic structures (Mlitz [21]), and they can also be
adopted for the non-algebraic structures. Furthermore, it is then shown that
congruences can be used to give a characterization of the connectednesses of
topological spaces and graphs which resembles the well-known and classical
characterizations of the algebraic radical classes ([32, 34].

\bigskip

\noindent The fact that connectednesses and disconnectednesses of
topological spaces and graphs can be defined as Hoehnke radicals, has opened
up a number of new questions to explore; sometimes with interesting
consequences. In the classical torsion theory, torsionfree classes (=
semisimple classes) are always hereditary and a hereditary torsion theory
means that the associated torsion class (= radical class) is hereditary. For
associative rings, the semisimple classes are also always hereditary, so
also here a hereditary radical would mean the associated radical class is
hereditary. But for more general classes, e.g. not necessarily associative
rings or near-rings, semisimple classes need not be hereditary. Thus, for a
radical in general, it is now customary to call it ideal-hereditary if both
its semisimple class and its radical class are hereditary. The relationships
between Hoehnke radicals, Kurosh-Amitsur radicals, connectednesses and
disconnectednesses and torsion theories with or without additional
properties (like hereditariness) have been investigated and clarified for
most concrete categories. For associative rings and similar types of
algebraic categories, in fact for $\Omega $-groups in general, it is
well-known that any ideal-hereditary Hoehnke radical is a Kurosh-Amitsur
radical. In terms of torsion theories, this statement says that any
hereditary torsion theory is Kurush-Amitsur. It is thus interesting to know
if an ideal-hereditary Hoehnke radical of topological spaces or graphs is a
Kurosh-Amistsur radical (i.e., whether every hereditary torsion theory of
topological spaces or graphs gives rise to a corresponding pair of
connectedness and disconnectedness). Wiegandt [38] has shown that for $%
\mathcal{S}$-acts this is not the case: a hereditary torsion theory need not
determine a Kurosh-Amitsur radical. It will be seen below that this is also
the case for topological spaces and graphs (see [31, 33]). In this vein of
the unexpected, we will find the radical theory of graphs that do not admit
loops. Also here there is a theory of congruences and Hoehnke radicals can
be defined. But they all degenerate, and this in spite of the fact that
there are non-trivial connecetednesses and disconnectednesses which, by the
way, all come as complementary pairs [36].

\bigskip

\noindent Our objecive here is to give an overview of the value added to the
general radical theory by this approach to the connectednesses and
disconnectednesses via congruences. To set the scene and to compare and
appreciate the correspondences and differences between the theories, the
appropriate definitions and main results on the general radical theory of
algebraic structures, actually mainly for associative rings, will be
recalled in the next section. In Sections 3, 4 and 5, one each for the
category of topological spaces, the category of graphs for which loops are
allowed and the category of graphs for which loops are not allowed
respectively, we present a brief summary of the congruence theory and the
radical theory for each, concluding with some interesting features of this
theory for that category. In this survey, the earlier works of the author
and his co-authors on the congruences of topological spaces and graphs and
their radical theory are freely used and often quoted verbatim without
reference. In particular, we refer to [4, 30, 31, 32, 33, 34, 36].

\section{Rings}

\noindent Here we give a quick overview of Kusorh-Amitsur radicals, Hoehnke
radicals, torsion theories and the relations betweeen them. The terminolgy
here is largely motivated by that in use for associative rings. Moreover,
the results recalled below are mostly for associative rings and may not
necessarily be valid for other classes or rings or algebras. Gardner and
Wiegandt [17] can be consulted for a thorough overview of the radical theory
of associative rings. To appreciate the similarities between the algebraic
and the non-algebraic with the use of congruences in the radical theory
presented in the next sections, [22] can consulted.

\noindent A property $\mathcal{R}$ that a ring may possess, is called a 
\textit{radical property} provided the following three conditions are
fulfilled:

$(R1)$ Any homomorphic image of a ring with property $\mathcal{R}$ has
property $\mathcal{R}.$

$(R2)$ Any ring has a largest ideal, called the $\mathcal{R}$\textit{%
-radical }of the ring, which, as a ring, has property

$\qquad \mathcal{R},$ and contains all other ideals of the ring which has
property $\mathcal{R}.$

$(R3)$ The quotient of any ring by its $\mathcal{R}$-radical contains no
non-trivial ideals with property $\mathcal{R}$

\qquad (such rings are called $\mathcal{R}$\textit{-semisimple}).

\noindent The class of all rings which has property $\mathcal{R}$ is usually
also denoted by $\mathcal{R}$ and no distinction is made between a ring in
the class $\mathcal{R}$ and a ring with the property $\mathcal{R}.$ If $%
\mathcal{R}$ is a radical property, then $\mathcal{R}$ is called a \textit{%
radical class}, the $\mathcal{R}$-radical of a ring $A$ is written as $%
\mathcal{R}(A)$ and the class $\mathcal{S}:=\{$ rings $A\mid \mathcal{R}%
(A)=0\}$ is called the \textit{semisimple class} of $\mathcal{R}.$ In
general, a class of rings is called a semisimple class if it is the
semisimple class of some radical class. In recognition of their
contributions in establishing this abstract approach to radicals by Kurosh
and independently Amitsur, these radicals are called $KA$\textit{-radicals}.
For a class of rings $\mathcal{M}$ and a ring $A,$ we let $\mathcal{M}%
(A):=\sum \{I\vartriangleleft A\mid I\in \mathcal{M}\}$ and $(A)\mathcal{M}:=%
\mathcal{\cap \{}I\vartriangleleft A\mathcal{\mid }A/I\mathcal{\in M\}.}$
The class $\mathcal{M}$ is said to be: \textit{hereditary} if for all rings $%
A,$ $I\lhd A\in \mathcal{M}$ implies $I\in \mathcal{M};$ \textit{%
homomorphically closed} if for any surjective homomorphism $\theta
:A\rightarrow B$ with $A\in \mathcal{M}$ implies $B\in \mathcal{M};$ \textit{%
inductive} if $I_{1}\subseteq I_{2}\subseteq ...$ is an ascending chain of
ideals of the ring $A$ with $I_{n}\in \mathcal{M}$ for all $n=1,2,3,...,$
then $\bigcup\limits_{n=1}^{\infty }I_{n}\in \mathcal{M;}$ \textit{closed
under extensions} if $I\lhd A$ with both $I$ and $A/I$ in $\mathcal{M},$
then also $A\in \mathcal{M;}$ \textit{regular} if every non-zero ideal of
the ring $A$ has a non-zero homomorphic image which is in $\mathcal{M},$
then $A\in \mathcal{M;}$ and \textit{closed under subdirect sums} whenever
rings $A_{j}\in \mathcal{M}$ for all $j\in J,$ then so is their subdirect
sum.

\begin{theorem}
For a class of rings $\mathcal{R},$ the four conditions $(1),(2),(3)$ and $%
(4)$ below are equivalent:

\noindent $(1)$ $\mathcal{R}$ is a KA-radical class.

\noindent $(2)$ $\mathcal{R}$ fulfills the following three conditions:

\qquad (a) $\mathcal{R}$ is homomorphically closed.

\qquad (b) $\mathcal{R}(A)\in \mathcal{R}$ for all rings $A.$

\qquad (c) $\mathcal{R}(A/\mathcal{R}(A))=0$ for all rings $A.$

\noindent $(3)$ $\mathcal{R}$ fulfills the following condition:

For any ring $A,$ $A\in \mathcal{R}\Leftrightarrow $ every non-zero
homomorphic image of $A$ has a non-zero ideal

which is in $\mathcal{R}.$

\noindent $(4)$ $\mathcal{R}$ fulfills the following three conditions:

\qquad (a) $\mathcal{R}$ is homomorphically closed.

\qquad (b) $\mathcal{R}$ is inductive.

\qquad (c) $\mathcal{R}$ is closed under extensions$.$
\end{theorem}

\noindent Semisimple classes can be characterized in their own right:

\begin{theorem}
For a class of rings $\mathcal{S},$ the four conditions $(1),(2),(3)$ and $%
(4)$ below are equivalent:

\noindent $(1)$ $\mathcal{S}$ is a KA-semisimple class.

\noindent $(2)$ $\mathcal{S}$ fulfills the following four conditions:

\qquad (a) $\mathcal{S}$ is regular.

\qquad (b) $\mathcal{S}$ is closed under subdirect sums$.$

\qquad (c) $\mathcal{S}$ is closed under extensions.

\qquad (d) $((A)\mathcal{S)S\rhd }A$ for all rings $A.$

\noindent $(3)$ $\mathcal{S}$ fulfills the following condition:

For any ring $A,$ $A\in \mathcal{S}\Leftrightarrow $ every non-zero ideal of 
$A$ has a non-zero homomorphic image

which is in $\mathcal{S}.$

\noindent $(4)$ $\mathcal{S}$ fulfills the following three conditions:

\qquad (a) $\mathcal{S}$ is hereditary.

\qquad (b) $\mathcal{S}$ is closed under subdirect sums.

\qquad (c) $\mathcal{S}$ is closed under extensions.
\end{theorem}

\noindent Condition $(3)$ in each of the previous two theorems has become
the defining condition of a radical class and a semisimple class
respectively in a wide-raging and diverse range of categories. Hoehnke's
approach to radicals [18] has prompted the following: an \textit{%
ideal-mapping} $\rho $ is a mapping which assigns to each ring $A$ an ideal $%
\rho (A),$ usually written as $\rho _{A}.$ The Hoehnke radical is defined as
follows:

\noindent An $H$\textit{-radical} is an ideal-mapping which fulfills the
following two conditions:

\qquad $(H1)$ For any surjective homomorphism $\theta :A\rightarrow B,$ $%
\theta (\rho (A))\subseteq \rho (B).$

\qquad $(H2)$ For all rings $A,$ $\rho (A/\rho (A))=0.$

\noindent For the $H$-radical $\rho ,$ $\mathcal{R}_{\rho }:=\{$rings $A\mid
\rho (A)=A\}$ is the associated radical class and $\mathcal{S}_{\rho }:=\{$%
rings $A\mid \rho (A)=0\}$ is the \textit{associated semisimple class}.

\begin{theorem}
Let $\mathcal{M}$ be a class of rings which contains the zero ring and is
closed under isomorphic copies. Then $\rho $ defined by $\rho (A)=\cap
\{I\rhd A\mid A/I\in \mathcal{M}\}$ for all rings $A$ is an $H$-radical with 
$\mathcal{S}_{\rho }$ the subdirect closure of the class $\mathcal{M}.$
Conversely, if $\rho $ is an $H$-radical, there is a class $\mathcal{M}$ of
rings which contains the zero ring, is closed under isomorphic copies, has $%
\rho (A)=\cap \{I\rhd A\mid A/I\in \mathcal{M}\}$ for all rings $A$ and $%
\mathcal{S}_{\rho }$ is the subdirect closure of $\mathcal{M}.$
\end{theorem}

\noindent To clarify the connection between KA-radicals and H-radicals, we
need two more definitions. An ideal mapping $\rho $ is: \textit{complete} if
for all rings $A,$ $\rho (I)=I\lhd A$ implies $I\subseteq \rho (A);$ and 
\textit{idempotent} if $\rho (\rho (A))=\rho (A)$ for all rings $A.$ In view
of the next result, a complete and idempotent $H$-radical is called a $KA$%
\textit{-radical}.

\begin{theorem}
For an ideal-mapping $\rho ,$ the following three conditions are equivalent:

\noindent (1) $\rho $ is an $H$-radical which is complete and idempotent.

\noindent (2) $\rho $ is an $H$-radical for which the corresponding radical
class $\mathcal{R}_{\rho }$ is a KA-radical class with associated semisimple
class $\mathcal{S}_{\rho }.$

\noindent (3) $\rho $ is an $H$-radical for which the corresponding
semisimple class $\mathcal{S}_{\rho }$ is a KA-semisimple class with
associated radical class $\mathcal{R}_{\rho }.$

\noindent Conversely, if $\mathcal{R}$ is a KA-radical class, then the ideal
mapping $\rho $ defined by $\rho (A)=\sum \{I\vartriangleleft A\mid I\in 
\mathcal{R}\}$ for all rings $A$ is a complete, idempotent H-radical with $%
\mathcal{R}_{\rho }=\mathcal{R;}$ or equivalently, if $\mathcal{S}$ is a
KA-semisimple class, then the mapping $\rho $ defined by $\rho (A)=\cap
\{I\rhd A\mid A/I\in \mathcal{S}\}$ for all rings $A$ is a complete,
idempotent H-radical with $\mathcal{S}_{\rho }=\mathcal{S.}$
\end{theorem}

\noindent An ordered pair $(\mathcal{R},\mathcal{S})$ of classes of rings $%
\mathcal{R}$ and $\mathcal{S}$ is called a \textit{torsion theory} if :

(T1) $\mathcal{R}\cap \mathcal{S}=0\mathcal{.}$

(T2) $\mathcal{R}$ is homomorphically closed.

(T3) $\mathcal{S}$ is hereditary.

(T4) Every ring $A$ has an ideal $B$ with $B\in \mathcal{R}$ and $A/B\in 
\mathcal{S.}$

\noindent The class $\mathcal{R}$ is called a \textit{torsion class} and $%
\mathcal{S}$ is a \textit{torsionfree class}. A \textit{hereditary torsion
theory} is a torsion theory $(\mathcal{R},\mathcal{S})$ with hereditary
torsion class $\mathcal{R}.$ This abstract approach to torsion theories was
introduced by Dickson [12] and has its origins in abelian categories.

\begin{theorem}
The following are equivalent:

\noindent $(1)$ $(\mathcal{R},\mathcal{S})$ is a torsion theory.

\noindent $(2)$ $\mathcal{R}$ is a KA-radical class with corresponding
semisimple class $\mathcal{S}$ and $\mathcal{R}(I)\subseteq \mathcal{R}(A)$
for all rings $A$ and $I\rhd A$.

\noindent $(3)$ $\mathcal{S}$ is a hereditary KA-semisimple class with
corresponding radical class $\mathcal{R}.$
\end{theorem}

\noindent An ideal mapping $\rho $ is called \textit{ideal-hereditary} if $%
\rho (I)=\rho (A)\cap I$ for all rings $A$ and $I\rhd A.$ Recall, $\rho $ is
a KA-radical if it is an idempotent and complete $H$-radical. For such a
radical $\rho ,$ let $\mathcal{R}$ and $\mathcal{S}$ be the corresponding
radical and semisimple class respectively. Then:

(1) $\mathcal{R}$ is hereditary iff $\rho (A)\cap I\subseteq \rho (I)$ for
all rings $A$ and $I\rhd A.$

(2) $\mathcal{S}$ is hereditary iff $\rho (I)\subseteq \rho (A)\cap I$ for
all rings $A$ and $I\rhd A.$

\noindent If $\rho $ is only a Hoehnke radical, then:

(3) $\rho (A)\cap I\subseteq \rho (I)$ for all rings $A$ and $I\rhd A$
implies $\rho $ is idempotent and $\mathcal{R}_{\rho }$ is hereditary.

(4) $\rho (I)\subseteq \rho (A)\cap I$ for all rings $A$ and $I\rhd A$
implies $\rho $ is complete and $\mathcal{S}_{\rho }$ is hereditary.

\noindent We know that an idempotent and complete H-radical of rings is a
KA-radical; hence an ideal-hereditary H-radical of associative rings
coincides with an ideal-hereditary KA-radical (and both the radical and
semisimple classes are hereditary) which in turn is the same as a hereditary
torsion theory.

\section{Topology}

\noindent In the early seventies, Preu\ss\ [23, 24, 25] defined and
developed a general theory of connectednesses and disconnectednesses for
topological spaces. This showed and clarified, amongst others, the
relationship between the separation axioms and non-connectedness on the one
hand and connected spaces on the other, establishing a Galois connection
between them. Subsequently many papers appeared on connectednesses and
disconnectednesses of topological spaces, see for example Castellini [6, 7,
8], Castellini and Holgate [9], Clementino [10] and Clementino and Tholen
[11] and their references. Of particular interest here are the results of
Arhangel'ski\u{\i} and Wiegandt [1]. They showed that the theory of
connectednesses and disconnectednesses of topological spaces is, from a
categorical perspective, the same as the KA-radical and -semisimple theory
of rings and related algebraic structures and also to the torsion theory of
abelian categories. In algebra, these radicals can be obtained as H-radicals
being the intersection of congruences (ideals) for which the quotients are
semisimple. In [30] it was shown that this can also be done for topological
spaces using the recently defined congruences for topological spaces: every
disconnectedness of topological spaces (and so also every connectedness) can
be obtained as an H-radical of the intersection of all congruences for which
the corresponding weak quotient is in the disconnectedness. Moreover, the
connectednesses of topological spaces can be characterized in terms of
congruences in precisely the same way that the radical classes of
associative rings are characterized in terms of ideals. But not all the
salient features of the radical theory of rings hold for their topological
counterparts. We shall see that a hereditary torsion theory (=
ideal-hereditary $H$-radical) of topological spaces need not be a
KA-radical. In fact, we give all the ideal-hereditary H-radicals; in
particular also showing which ones are KA-radicals. Details of the work
presented here can be found in [30, 31, 32].

\bigskip

\noindent We start with an overview of the congruence theory for topological
spaces. A topological space will usually be denoted by $(X,\mathcal{T})$ and
when there is no need to specify the topology, just by $X.$ The one-element
space will be denoted by $T$ and we identify all one-element spaces with $T.$
A \textit{trivial space }is either $T$ or the empty space $\emptyset .$ A
subset of a topological space will always be regarded as a topological space
with respect to the relative topology, unless explicitly mentioned
otherwise. Homeomorphic topological spaces may be denoted by $\cong .$

\subsection{\noindent Topological congruences}

\begin{definition}
Let $(X,\mathcal{T})$ be a topological space. A congruence $\rho $ on $X$ is
a pair $\rho =(\sim ,\mathbb{T})$ where:

\noindent (C1) $\sim $ is an equivalence relation on $X.$

\noindent (C2) $\mathbb{T}$ is a topology on $X$ with $\mathbb{T}\subseteq 
\mathcal{T.}$ $\mathbb{T}$ is called the congruence topology.

\noindent (C3) For all $x\in X,$ $x\in U\in \mathbb{T}$ implies $%
[x]\subseteq U.$

\noindent A congruence $\rho =(\sim ,\mathbb{T})$ is a strong congruence on $%
X$ if $\mathbb{T=\{}U\subseteq X\mathbb{\mid }U$ is open in $X$ and $x\in U$
implies $[x]\subseteq U\}.$
\end{definition}

\noindent Here $[x]$ denotes the equivalence class of $x\in X$ with respect
to the equivalence relation $\sim $ on $X.$ By $(C3),$ if $U\in \mathbb{T}$,
then $U=\bigcup\limits_{a\in U}[a]$ and if the congruence is strong, then $%
\mathbb{T}$ consists precisely of all open sets $U$ with this property. As
examples of congruences on $(X,\mathcal{T})$ we mention the \textit{identity
congruence} on $X,$ $\iota _{X}=(\Bumpeq ,\mathcal{T})$ where $a\Bumpeq b$
if and only if $a=b,$ the \textit{universal congruence} on $X,$ $\upsilon
_{X}=(\leftrightsquigarrow ,\mathcal{I}_{X})$ where $\leftrightsquigarrow $
is the universal equivalence relation on $X$ with $a\leftrightsquigarrow b$
for all $a,b\in X$ and $\mathcal{I}_{X}$ is the indiscrete topology on $X$
and the \textit{kernel }of a continuous map $f:X\rightarrow Y,$ denoted by $%
\ker f=(\sim _{f},\mathbb{T}_{f}),$ with $a\sim _{f}b\Leftrightarrow $ $%
f(a)=f(b)$ for $a,b\in X$ and $\mathbb{T}_{f}=\{f^{-1}(V)\mid V\subseteq Y$
open\}. For this congruence $[x]=f^{-1}(f(x))$ for all $x\in X$ and for any $%
U\in \mathbb{T}_{f},$ $f^{-1}(f(U))=U.$ The \textit{strong kernel of }$f,$
denoted by sker $f=(\sim _{f},\mathbb{T}_{sf}),$ is the congruence on $X$
with $\sim _{f}$ as defined above and $\mathbb{T}_{sf}=\{U\subseteq X\mid U$
is open in $X$ and $U=f^{-1}(f(U))\}.$ Note that for a topological space $(X,%
\mathcal{T})$ and any topology $\mathbb{T}$ on $X$ with $\mathbb{T}\subseteq 
\mathcal{T,}$ $(\Bumpeq ,\mathbb{T)}$ is a congruence on $X$ and is called a 
\textit{trivial congruence} on $X.$ A trivial congruence on $X$ need not be
the identity congruence $\iota _{X},$ but a strong congruence is trivial if
and only if it coincides with the identity congruence $\iota _{X}.$ Every
congruence $\rho =(\sim ,\mathbb{T})$ on $X$ determines a topological space $%
(X/\rho ,\mathcal{T}/\rho )$ with $X/\rho =\{[x]\mid x\in X\}$ and topology $%
\mathcal{T}/\rho =\{\pi _{\rho }(U)\mid U\in \mathbb{T\}}$ where $\pi _{\rho
}:X\rightarrow X/\rho $ defined by $\pi _{\rho }(x)=[x]$ is a surjective
continuous map with $\ker \pi _{\rho }=\rho .$ This space $X/\rho $ is
called the \textit{weak quotient space determined by }$\rho $ and $\pi
_{\rho }$ is the weak quotient map (or just called the canonical map). In
general $X/\rho $ need not be a quotient space with $\pi _{\rho }$ a
quotient map. But if $\rho $ is a strong congruence, $X/\rho $ is a quotient
space and $\pi _{\rho }$ is a quotient map. Using condition $(C3)$ and the
surjectivity of the map $\pi _{\rho },$ it can be shown that $\mathcal{T}%
/\rho =\{W\subseteq X/\rho \mid \pi _{\rho }^{-1}(W)\in \mathbb{T\}.}$ Two
expected quotients are $(X/\iota _{X},\mathcal{T}/\iota _{X})\cong (X,%
\mathcal{T})$ and for $X\neq \emptyset ,$ $(X/\upsilon _{X},\mathcal{T}%
/\upsilon _{X})\cong T.$ In fact, it can be shown that $X/\rho \cong
X\Leftrightarrow \rho =\iota _{X}$; for $X\neq \emptyset ,$ $X/\rho \cong
T\Leftrightarrow \rho =\upsilon _{X};$ and for any $X,$ $\iota _{X}=\upsilon
_{X}\Leftrightarrow X\cong T$ or $X=\emptyset .$

\bigskip

\noindent \textbf{Ordering of congruences. }For two congruences $\rho =(\sim
_{\rho },\mathbb{T}_{\rho })$ and $\gamma =(\sim _{\gamma },\mathbb{T}%
_{\gamma })$ on $X,$ $\rho $\textit{\ is contained in} $\gamma ,$ written as 
$\rho \sqsubseteq \gamma ,$ provided $\sim _{\rho }$ $\subseteq $ $\sim
_{\gamma }$ and $\mathbb{T}_{\gamma }\subseteq \mathbb{T}_{\rho }.$ For any
congruence $\rho $ on $X,$ $\iota _{X}\sqsubseteq \rho \sqsubseteq \upsilon
_{X}$ and $\rho \sqsubseteq \rho .$ This ordering $\sqsubseteq $ is a
partial order on the class $Con(X):=\{\theta \mid \theta $ is a congruence
on $X\}.$ In fact, $Con(X)$ is a bounded complete lattice. The meet of the
congruences $\theta _{i}=(\sim _{i},\mathbb{T}_{i})\in Con(X),$ $i\in I$, is
given by the intersection of congruences $\bigcap\limits_{i\in I}\theta _{i}$
$=(\sim _{\cap },\mathbb{T}_{\cap })$ defined by: For $a,b\in X,$ $a\sim
_{\cap }b\Leftrightarrow a\sim _{i}b$ for all $i\in I$ and the congruence
topology $\mathbb{T}_{\cap }$ is given by the topology on $X$ with basis $%
\mathcal{B}:=\{B\subseteq X\mid B$ is a finite intersection $%
B=\bigcap\limits_{j=1}^{n}U_{j}$ where $U_{j}\in \mathbb{T}_{i_{j}}$ for
some $i_{j}\in I,j=1,2,3,...,n,n\geq 1$\}. The join of the congruences $%
\theta _{i}=(\sim _{i},\mathbb{T}_{i})$ is given by the sum $%
\sum\limits_{i\in I}\theta _{i}=(\sim _{\Sigma },\mathbb{T}_{\Sigma })$
defined by: $\mathbb{T}_{\Sigma }:=\bigcap\limits_{i\in I}\mathbb{T}_{i}$
and for $a,b\in X,$

$a\sim _{\Sigma }b\Leftrightarrow $ there are $i_{1},i_{2},...,i_{n}\in I$
and $a_{i_{1}},a_{i_{2}},...,a_{i_{n}}\in X,n\geq 2,$ such that $%
a=a_{i_{1}}\sim _{i_{1}}a_{i_{2}}\sim _{i_{2}}a_{i_{3}}\sim _{i_{3}}...\sim
_{i_{n-2}}a_{i_{n-1}}\sim _{i_{n-1}}a_{i_{n}}=b.$

\noindent Note that if $\circ $ denotes the usual composition of two binary
relations, then

\begin{equation*}
a\sim _{\Sigma }b\Leftrightarrow \ a(\sim _{i_{1}}\circ \sim _{i_{2}}\circ
\sim _{i_{3}}\circ ...\circ \sim _{i_{n-1}}\circ \sim _{i_{n}})b\text{ for
some }n\geq 1,i_{1},i_{2},...,i_{n}\in I.
\end{equation*}%
\bigskip

\noindent If $\theta _{i}$ is a strong congruence for all $i,$ then it can
easily be shown that the sum $\sum\limits_{i\in I}\theta _{i}$ is also a
strong congruence. The sum of two congruences $\alpha $ and $\beta $ on a
space $X$ is written as $\alpha +\beta .$ The usual relationships are valid: 
$\alpha +\iota _{X}=\alpha =\alpha \cap \upsilon _{X},$ $\alpha \cap \iota
_{X}=\iota _{X},\alpha +\upsilon _{X}=\upsilon _{X}$ and $\alpha \sqsubseteq
\beta \Leftrightarrow \alpha =\alpha \cap \beta \Leftrightarrow \beta
=\alpha +\beta .$ Since $x\sim x$ for any equivalence $\sim ,$ there is no
loss of generality in writing $a\sim _{\alpha +\beta }b$ as $a=a_{i_{1}}\sim
_{\alpha }a_{i_{2}}\sim _{\beta }a_{i_{3}}\sim _{\alpha }...\sim _{\alpha
}a_{i_{n-1}}\sim _{\beta }a_{i_{n}}=b$ for some $n\geq 1$ and $%
a_{i_{1}},a_{i_{2}},...,a_{i_{n}}\in X.$

\bigskip

\noindent \textbf{Homeomorphism Theorems.}

\begin{theorem}
(First Homeomorphism Theorem) Let $f:(X,\mathcal{T})\rightarrow (Y,\mathcal{F%
})$ be a surjective continuous map with $\alpha =\ker f.$ Then $(X/\alpha ,%
\mathcal{T}/\alpha )\cong (Y,\mathcal{F}).$
\end{theorem}

\noindent Let $(X,\mathcal{T})$ be a topological space with $\rho =(\sim ,%
\mathbb{T})$ a congruence on $X$ and $\pi _{\rho }:X\rightarrow X/\rho $ the
canonical map $\pi _{\rho }(a)=[a]$. Let $S$ be a non-empty subset of $X.$
Then $\rho $ induces a congruence on the subspace $S$, denoted by $S\cap
\rho =(\sim _{S},\mathbb{T}_{S}),$ by restricting $\rho $ to $S$ in a
natural way: For all $a,b\in S,$ $a\sim _{S}b\Leftrightarrow a\sim b$ and $%
\mathbb{T}_{S}=\{U\cap S\mid U\in \mathbb{T\}}.$ This congruence $S\cap \rho 
$ is called the \textit{restriction of the congruence }$\rho $\textit{\ to }$%
S.$ The subspace of $X/\rho $ determined by $\pi _{\rho }(S)$ will be
denoted by $(S+\rho )/\rho $.

\begin{theorem}
(Second Homeomorphism Theorem) Let $(X,\mathcal{T})$ be a topological space, 
$S$ a subspace of $X$ and $\rho =(\sim ,\mathbb{T})$ a congruence on $X.$
Then the weak quotient space of $S$ determined by the congruence $S\cap \rho 
$ on $S$ is homeomorphic to the subspace $\{[a]\mid a\in S\}$ of $(X/\rho ,%
\mathcal{T}/\rho );$ i.e., $S/S\cap \rho \cong S+\rho /\rho .$
\end{theorem}

\noindent For the next homeomorphism theorem, we need the quotient of two
congruences. Let $(X,\mathcal{T})$ be a topological space; let $\alpha
=(\sim _{\alpha },\mathbb{T}_{\alpha })$ and $\beta =(\sim _{\beta },\mathbb{%
T}_{\beta })$ be two congruences on $X$ with $\alpha \sqsubseteq \beta .$
The \textit{quotient of }$\beta $\textit{\ by }$\alpha $, written as $\beta
/\alpha =(\sim _{\beta /\alpha },\mathbb{T}_{\beta /\alpha }),$ is the
congruence on the weak quotient space $(X/\alpha ,\mathcal{T}/\alpha )$
defined as follows: for $[a]_{\alpha },[b]_{\alpha }\in X/\alpha $, $%
[a]_{\alpha }\sim _{\beta /\alpha }[b]_{\alpha }\Leftrightarrow a\sim
_{\beta }b$ and $\mathbb{T}_{\beta /\alpha }=\{W\subseteq X/\alpha \mid
W=\pi _{\alpha }(U)$ for some $U\in \mathbb{T}_{\beta }\}$ where $\pi
_{\alpha }:X\rightarrow X/\alpha $ is the canonical map. In fact we have:

\begin{lemma}
Let $(X,\mathcal{T})$ be a topological space with $\alpha =(\sim _{\alpha },%
\mathbb{T}_{\alpha })$ a congruence on $X.$ Then $\gamma =(\sim ,\mathbb{T})$
is a congruence on the weak quotient space $(X/\alpha ,\mathcal{T}/\alpha )$
if and only if $\gamma =\beta /\alpha $ for some congruence $\beta =(\sim
_{\beta },\mathbb{T}_{\beta })$ on $X$ with $\alpha \sqsubseteq \beta .$
\end{lemma}

\begin{theorem}
(Third Homeomorphism Theorem) Let $(X,\mathcal{T})$ be a topological space
and let $\alpha =(\sim _{\alpha },\mathbb{T}_{\alpha })$ and $\beta =(\sim
_{\beta },\mathbb{T}_{\beta })$ be two congruences on $X$ with $\alpha
\sqsubseteq \beta .$ Then $\beta /\alpha $ is a congruence on $X/\alpha $
and $(X/\alpha )/(\beta /\alpha )$ is homeomorphic to $X/\beta .$
\end{theorem}

\begin{corollary}
Let $\theta $ be a congruence on $X$. Then there is a one-to-one
correspondence between the set of all congruences $\alpha $ on $X$ for which 
$\theta \sqsubseteq \alpha $ and the set of all congruences on $X/\theta $
given by $\alpha \rightarrow \alpha /\theta .$ This correspondence preserves
containment, joins and intersections.
\end{corollary}

\bigskip

\noindent \textbf{Image of a congruence. }\noindent Let $f:X\rightarrow Y$
be a surjective continuous map with $\rho =(\sim _{\rho },\mathbb{T}_{\rho
}) $ a congruence on $X$. Then $f(\rho )=(\sim _{f(\rho )},\mathbb{T}%
_{f(\rho )})$ is the congruence on $Y=f(X)$ defined by $f(a)\sim _{f(\rho
)}f(b)\Leftrightarrow $ there are $a_{1},a_{2},...,a_{n}$ in $X$ with $a\sim
_{\rho }a_{1},$ $f(a_{1})=f(a_{2}),$ $a_{2}\sim _{\rho }a_{3},$ $%
f(a_{3})=f(a_{4}),$ $a_{4}\sim _{\rho }a_{5},....,a_{n-1}\sim _{\rho }a_{n},$
$f(a_{n})=f(b)$ and $\mathbb{T}_{f(\rho )}=\{V\subseteq Y$ open$\mid $ $%
f^{-1}(V)\in \mathbb{T}_{\rho }\}.$ It is clear that $f(a)\sim _{f(\rho
)}f(b)\Leftrightarrow a\sim _{\rho +\ker f}b\Leftrightarrow a^{\prime }\sim
_{\rho +\ker f}b^{\prime }$ for $f(a^{\prime })=f(a)$ and $f(b^{\prime
})=f(b).$ We note that if $\rho $ is strong, so is $f(\rho )$ and in view of
the First Homeomorphism Theorem, if $\alpha =\ker f,$ then $f(\rho )=(\rho
+\alpha )/\alpha $ where $Y\cong X/\alpha $. As is to be expected, $f(\iota
_{X})=\iota _{Y}$ and $f(\upsilon _{X})=\upsilon _{Y}.$

\bigskip

\noindent \textbf{Subdirect products. }Let $\prod\limits_{i\in I}X_{i}$
denote the product of the topological spaces $X_{i},i\in I$ with $p_{j}:$ $%
\prod\limits_{i\in I}X_{i}\rightarrow X_{j}$ the $j$-th projection. A
subspace $Y$ of $\prod\limits_{i\in I}X_{i}$ is said to be a \textit{%
subdirect product} of the spaces $X_{i},i\in I,$ if $p_{i}(Y)=X_{i}$ for all 
$i\in I.$ As is the case for algebra, subdirect products can be
characterized in terms of congruences:

\begin{theorem}
A topological space $Y$ is a subdirect product of spaces $X_{i},i\in I,$ if
and only if for every $i\in I$ there is a congruence $\theta _{i}$ on $Y$
such that $Y/\theta _{i}\cong X_{i}$ and $\bigcap\limits_{i\in I}\theta
_{i}=\iota _{Y}$.
\end{theorem}

\subsection{\noindent Radical theory}

\noindent As is usual in radical theory, all considerations will be in a 
\textit{universal class} $\mathcal{W}$ of topological spaces. This means $%
\mathcal{W}$ is a non-empty class of spaces which is \textit{hereditary} (if 
$Y$ is a subspace of $X\in \mathcal{W},$ then $Y\in \mathcal{W})$ and 
\textit{closed under continuous images} (if $f:X\rightarrow Y$ is a
surjective continuous map with $X\in \mathcal{W},$ then also $Y\in \mathcal{W%
}).$ Clearly then, $\mathcal{W}$ contains the trivial spaces (empty space
and one-point spaces). By assumption, any subclass of $\mathcal{W}$ under
discussion, will be an \textit{abstract class}; i.e., it contains the
trivial spaces and all homeomorphic copies of spaces from the class (this
assumption will mostly not even be mentioned explicitly).

\bigskip

\noindent From [23], we recall: a class $\mathcal{C}$ of topological spaces
in $\mathcal{W}$ is called a \textit{connectedness }if there is a class of
spaces $\mathcal{P}$ in $\mathcal{W}$ such that $\mathcal{C}=\{X\in \mathcal{%
W}\mid $every continuous mapping $X\rightarrow Y\in \mathcal{P}$ is
constant\} and a class of spaces $\mathcal{D}$ in $\mathcal{W}$ is called a 
\textit{disconnectedness }if there is a class of spaces $\mathcal{Q}$ in $%
\mathcal{W}$ such that $\mathcal{D}=\{X\in \mathcal{W}\mid $every continuous
mapping $Y\rightarrow X$ with $Y\in \mathcal{Q}$ is constant\}.

\bigskip

\noindent Let $\mathcal{C}$ and $\mathcal{D}$ be be subclasses of $\mathcal{%
W.}$ A subspace $Y$ of a space $X$ is called a $\mathcal{C}$\textit{%
-subspace }of $X$ if $Y\in \mathcal{C.}$ A continuous image $Y$ of a space $%
X $ is called a \textit{non-trivial continuous image} of $X$ provided $Y$ is
not a trivial space. Likewise, a subspace $S$ of $X$ is a \textit{%
non-trivial subspace} of $X$ if it is not a trivial space. Motivated by the
terminology from ring theory, a class $\mathcal{C}$ of spaces is called a%
\textit{\ KA-radical class} if it satisfies: $X\in \mathcal{C}%
\Leftrightarrow $ every non-trivial continuous image of $X$ has a
non-trivial $\mathcal{C}$-subspace. $\mathcal{D}$ is a\textit{\
KA-semisimple class} if it satisfies: $X\in \mathcal{D}\Leftrightarrow $
every non-trivial subspace of $X$ has a non-trivial continuous image which
is in $\mathcal{D}.$ Arhangel'ski\u{\i} and Wiegandt [1] have shown that a
class of spaces $\mathcal{C}$ is a connectedness if and only if it is a
KA-radical class and a class of spaces $\mathcal{D}$ is a disconnectedness
if and only if it is a KA-semisimple class. We explicitly recall these two
statements for later reference.

\begin{proposition}
\lbrack 1] Let $\mathcal{C}$ and $\mathcal{D}$ be abstract classes of
topological spaces in $\mathcal{W}$.

\noindent (1) $\mathcal{C}$ is a connectedness if and only if $\mathcal{C}$
satisfies the condition: $X\in \mathcal{C}$ if and only if every non-trivial
continuous image of $X$ has a non-trivial $\mathcal{C}$-subspace.

\noindent (2) $\mathcal{D}$ is a disconnectedness if and only if $\mathcal{D}
$ satisfies the condition: $X\in \mathcal{D}$ if and only if every
non-trivial subspace of $X$ has a non-trivial continuous image in $\mathcal{D%
}$.
\end{proposition}

\noindent Two operators $\mathcal{U}$ and $\mathcal{S}$ on a class $\mathcal{%
M}\subseteq \mathcal{W},$ called the \textit{upper radical operator} and 
\textit{semisimple operator} respectively, are defined by: $\mathcal{UM}%
=\{X\in \mathcal{W}\mid X$ has no non-trivial continuous image in $\mathcal{M%
}\}$ and $\mathcal{SM}=\{X\in \mathcal{W}\mid X$ has no non-trivial subspace
in $\mathcal{M}\}.$ Note that $\mathcal{M}\cap \mathcal{UM}=\{T,\emptyset \}=%
\mathcal{M}\cap \mathcal{SM}$. The class $\mathcal{UM}$ is always closed
under continuous images and $\mathcal{SM}$ is always hereditary. If $%
\mathcal{M}$ is hereditary, then $\mathcal{UM}$ is a connectedness and if $%
\mathcal{M}$ is closed under continuous images, then $\mathcal{SM}$ is a
disconnectedness. It can be shown that any connectedness $\mathcal{C}$ is
closed under continuous images; hence $\mathcal{SC}$ is a disconnectedness,
called the \textit{disconnectedness} \textit{corresponding to} $\mathcal{C}.$
Likewise, any disconnectedness $\mathcal{D}$ is hereditary, hence $\mathcal{%
UD}$ is a connectedness called the \textit{connectedness corresponding to} $%
\mathcal{D}$. Moreover, $\mathcal{C}\subseteq \mathcal{W}$ is a\textit{\
connectedness }if and only if $\mathcal{C=USC}$ and $\mathcal{D}\subseteq 
\mathcal{W}$ is a\textit{\ }diconnectedness\textit{\ }if and only if $%
\mathcal{D}=\mathcal{SUD}$.

\bigskip

\noindent A mapping $\sigma $ which assigns to each $X\in \mathcal{W}$ a
congruence $\sigma (X)=\sigma _{X}=(\sim _{\sigma _{X}},\mathbb{T}_{\sigma
_{X}})$ on $X,$ is called a \textit{H-radical }on $\mathcal{W}$ if it
satisfies the following two conditions:

\qquad (H1) For any surjective continuous map $f:X\rightarrow Y,$ $f(\sigma
_{X})\sqsubseteq \sigma _{Y}.$

\qquad (H2) For all $X\in \mathcal{W}$, $\sigma (X/\sigma _{X})=\iota
_{X/\sigma _{X}},$ the identity congruence on $X/\sigma _{X}.$

\noindent The class $\mathcal{S}_{\sigma }=\{X\in \mathcal{W}\mid \sigma
_{X}=\iota _{X}\}$ is called the associated semisimple class and $\mathcal{R}%
_{\sigma }=\{X\in \mathcal{W}\mid \sigma _{X}=\upsilon _{X}\}$ the
associated radical class. Note that $\mathcal{S}_{\sigma }\cap \mathcal{R}%
_{\sigma }=\{T,\emptyset \},$ $\mathcal{R}_{\sigma }$ is always closed under
continuous images and $\mathcal{R}_{\sigma }=\mathcal{U}\mathcal{S}_{\sigma
}.$ Hoehnke radicals are very general as is shown in the next result which
also gives most of the salient features of these radicals:

\begin{theorem}
(1) Let $\sigma $ be an H-radical on $\mathcal{W}.$ Then, for every $X\in 
\mathcal{W},$ $\sigma (X)=\cap \{\theta \mid \theta $ is a congruence on $X$
for which $X/\theta \in \mathcal{S}_{\sigma }\}$ and $\mathcal{S}_{\sigma }$
is closed under subdirect products.

\noindent (2) Conversely, let $\mathcal{M}\subseteq \mathcal{W}$ be any
abstract class. Then $\sigma $ defined by $\sigma (X)=\cap \{\theta \mid
\theta $ is a congruence on $X$ for which $X/\theta \in \mathcal{M}\}$ for
all $X\in \mathcal{W}$ is an H-radical on $\mathcal{W}$ and $\mathcal{S}%
_{\sigma }=\overline{\mathcal{M}}$, the subdirect closure of $\mathcal{M}.$
\end{theorem}

\noindent In [30] it was shown that any disconnectedness (and hence also
every connectedness) can be obtained from an H-radical provided it fulfills
two additional requirements. An H-radical $\sigma $ on $\mathcal{W}$ is:

\noindent \textit{complete} if $X\in \mathcal{W}$ and $\theta $ is a strong
congruence on $X$ with $[a]_{\theta }\in \mathcal{R}_{\sigma }$ for all $%
a\in X,$ then $\theta \sqsubseteq \sigma _{X};$ and

\noindent \textit{idempotent} if for all $X\in \mathcal{W}$ and $a\in X,$ we
have $[a]_{\sigma _{X}}\in \mathcal{R}_{\sigma }.$

\noindent Then:

\begin{theorem}
(1) Let $\sigma $ be an H-radical on $\mathcal{W}$. Suppose $\sigma $ is
complete, idempotent and for all $X\in \mathcal{W}$, the congruence $\sigma
_{X}$ is strong. Then $\mathcal{S}_{\sigma }$ is a disconnectedness and $%
\mathcal{R}_{\sigma }=\mathcal{U}\mathcal{S}_{\sigma }$ is a connectedness.

\noindent (2) Let $\mathcal{D}\subseteq \mathcal{W}$ be a disconnectedness.
Then there is an H-radical $\sigma $ on $\mathcal{W}$ which is complete,
idempotent, for all $X\in \mathcal{W}$ the congruence $\sigma _{X}$ is
strong, $\mathcal{S}_{\sigma }=\mathcal{D}$ and $\mathcal{R}_{\sigma }=%
\mathcal{UD}$. The congruence $\sigma _{X}$ is given by $\sigma _{X}=\cap
\{\theta \mid \theta $ is a congruence on $X$ for which $X/\theta \in 
\mathcal{D}\}.$
\end{theorem}

\noindent In view of this result, a Hoehnke radical $\sigma $ on $\mathcal{W}
$ which is complete, idempotent and for which $\sigma _{X}$ is a strong
congruence for all $X$ is called a \textit{KA-radical. N}ext we give a
characterization of a connectedness in terms of congruences.We start by
recalling two fundamental theorems for connectednesses and
disconnectednesses from [1]. For this we firstly fix some terminology. Let $%
\mathcal{C}$ be a class of topological spaces. Then $\mathcal{C}$ is \textit{%
second additive} if, whenever a space $X$ is covered by a family $\mathcal{F}
$ of $\mathcal{C}$-subspaces with $\cap \mathcal{F}\neq \emptyset ,$ then $%
X\in \mathcal{C};$ and $\mathcal{C}$ is\textit{\ }$q$\textit{-reversible} if 
$f:X\rightarrow Y$ is a quotient map with $Y$ and $f^{-1}(y)$ in $\mathcal{C}
$ for all $y\in Y,$ then also $X\in \mathcal{C.}$ Recall, a subspace $Y$ of
a space $X$ is a $\mathcal{C}$\textit{-subspace }if $Y\in \mathcal{C}.$

\begin{theorem}
\lbrack 1] Let $\mathcal{C}$ be an abstract class of topological spaces in $%
\mathcal{W}$. Then statements $(1)$ and $(2)$ below are equivalent$:$

\noindent (1) $\mathcal{C}$ is a connectedness.

\noindent (2) $\mathcal{C}$ satisfies the following three conditions:

(a) $\mathcal{C}$ is closed under continuous images.

(b) $\mathcal{C}$ is second additive.

(c) $\mathcal{C}$ is $q$-reversible.
\end{theorem}

\begin{theorem}
\lbrack 1] Let $\mathcal{D}\subseteq \mathcal{W}$ be a\textit{\
disconnectedness }with associated connectedness $\mathcal{C}=\mathcal{UD}$.
For any $X\in \mathcal{W}$:

\noindent (1) There is a quotient map $q_{X}:X$ $\rightarrow $ $X_{s}$ with $%
X_{s}\in \mathcal{D}.$

\noindent (2) For every surjective continuous map $g:X\rightarrow Y\in 
\mathcal{D}$, there is a continuous map $h:X_{s}\rightarrow Y$ such that $%
h\circ q_{X}=g$.

\noindent (3) For every $t\in X_{s},$ $q_{X}^{-1}(t)\in \mathcal{C}.$

\noindent (4) If $Y$ is a subspace of $X$ with $Y\in \mathcal{C}$, then $%
Y\subseteq q_{X}^{-1}(t)$ for some $t\in X_{s}.$
\end{theorem}

\noindent For an abstract class $\mathcal{C}$ of topological spaces, a
congruence $\rho $ on $X$ is a $\mathcal{C}$\textit{-congruence} if it is a
strong congruence and for each $x\in X,$ $[x]\in \mathcal{C}$ holds. It will
be useful to rephrase Theorem 3.12 above in terms of congruences.

\begin{theorem}
\lbrack 1] Let $\mathcal{D}\subseteq \mathcal{W}$ be a\textit{\
disconnectedness} with associated connectedness $\mathcal{C}=\mathcal{UD.}$
For any $X\in \mathcal{W}$:

\noindent (1) The congruence $\sigma _{X}=\cap \{\theta \mid \theta $ is a
congruence on $X$ for which $X/\theta \in \mathcal{D}\}$ is a strong
congruence and the weak quotient map $\pi _{X}:X$ $\rightarrow $ $X/\sigma
_{X}$ is a quotient map with $X/\sigma _{X}\in \mathcal{D}.$

\noindent (2) For every surjective continuous map $g:X\rightarrow Y\in 
\mathcal{D}$, there is a continuous map $h:X/\sigma _{X}\rightarrow Y$ such
that $h\circ \pi _{X}=g$ (or equivalently, $\sigma _{X}\sqsubseteq \ker g).$

\noindent (3) $\sigma _{X}$ is a $\mathcal{C}$-congruence on $X.$

\noindent (4) If $Y$ is a $\mathcal{C}$-subspace of $X$, then $Y\subseteq
\lbrack x]_{\sigma _{X}}$ for some $x\in X.$
\end{theorem}

\noindent For our characterization of connectednesses, we start with:

\begin{proposition}
Let $\mathcal{C}$ an be abstract class of topological spaces in $\mathcal{W}%
. $ Then $\mathcal{C}$ is a connectedness if and only if $\mathcal{C}$
satisfies the condition: $X\in \mathcal{C}$ if and only if every non-trivial
continuous image of $X$ has a non-trivial $\mathcal{C}$-congruence.
\end{proposition}

\noindent Let $\mathcal{C}$ be an abstract class of topological spaces. For
any space $X,$ let $\rho (X)$ be the congruence on $X$ defined by $\rho
(X)=\sum \{\alpha \mid \alpha $ is a $\mathcal{C}$-congruence on $X\}$.
Mostly we will write $\rho (X)$ as $\rho _{X}$ but there are occasions when
the former will be better to use. A sum of strong congruences is a strong
congruence, hence $\rho _{X}$ is a strong congruence and it contains all $%
\mathcal{C}$-congruences of $X.$ Note that $\upsilon _{X}$ is a $\mathcal{C}$%
-congruence on $X$ if and only if $X\in \mathcal{C}.$ When $\mathcal{C}$ is
a connectedness, then more can be said as will be seen below.

\begin{proposition}
Let $\mathcal{C}$ be a connectedness with associated disconnectedness $%
\mathcal{D}=\mathcal{SC}$ in $\mathcal{W.}$ Then $\rho _{X}$ is a $\mathcal{C%
}$-congruence and $\rho _{X}=\sigma _{X}$ for all $X$ where $\rho _{X}=\sum
\{\alpha \mid \alpha $ is a $\mathcal{C}$-congruence on $X\}$ and $\sigma
_{X}=\cap \{\theta \mid \theta $ is a congruence on $X$ with $X/\theta \in 
\mathcal{D}\}.$
\end{proposition}

\noindent A class $\mathcal{C}$ of topological spaces is \textit{closed
under extensions} if it satisfies: whenever $\alpha $ and $\beta $ are
congruences on $X$ with $\alpha \sqsubseteq \beta ,$ $\alpha $ a $\mathcal{C}
$-congruence on $X$ and $\beta /\alpha $ a $\mathcal{C}$-congruence on $%
X/\alpha ,$ then $\beta $ is a $\mathcal{C}$-congruence on $X.$ This brings
us to the main result of this section which shows that a connectedness of
topological spaces can be characterized in terms of the congruence $\rho
_{X} $ and also in terms of $\mathcal{C}$-congruences (for associative
rings, one would say in terms of the radical and in terms of the radical
ideals respectively). In particular, this shows that the characterization of
connectednesses of topological spaces is in complete harmony with the
characterization of the radical classes of associative rings (compare the
next theorem with Theorem 2.1 above). A last concept we need is: A class $%
\mathcal{C}$ of topological spaces has the \textit{inductive property} if it
satisfies the following condition: whenever $\alpha _{1}\sqsubseteq \alpha
_{2}\sqsubseteq \alpha _{3}\sqsubseteq ...$ is a chain of $\mathcal{C}$%
-congruences on a topological space $X\in \mathcal{W}$, then $\sum \{\alpha
_{i}\mid i=1,2,3,...$ $\}$ is a $\mathcal{C}$-congruence on $X.$

\begin{theorem}
Let $\mathcal{C}$ be an abstract class of spaces in $\mathcal{W}$. Then
statements (1),(2) and (3) below are equivalent$:$

\noindent (1) $\mathcal{C}$ is a connectedness.

\noindent (2) $\mathcal{C}$ satisfies the following three conditions:

(a) For every surjective continuous map $f:X\rightarrow Y,$ $f(\rho
_{X})\sqsubseteq \rho _{Y}.$

(b) For every space $X,$ $\rho _{X}$ is a $\mathcal{C}$-congruence on $X$
and it contains all $\mathcal{C}$-congruences on $X.$

(c) For every space $X,$ $\rho (X/\rho _{_{X}})=\iota _{X/\rho _{X}}.$

\noindent (3) $\mathcal{C}$ satisfies the following three conditions:

(a) If $\alpha $ and $\beta $ are congruences on $X$ with $\alpha
\sqsubseteq \beta $ and $\beta $ a $\mathcal{C}$-congruence$,$ then $\beta
/\alpha $ is a $\mathcal{C}$-con-

\qquad gruence on $X/\alpha .$

(b) $\mathcal{C}$ has the inductive property.

(c) $\mathcal{C}$ is closed under extensions.
\end{theorem}

\subsection{Hereditary torsion theories.}

\noindent For this section, the universal class $\mathcal{W}$ is the class
of all topological spaces. For any set $X,$ we will use $\mathcal{I}_{X}$
and $\mathcal{D}_{X}$ to denote the indiscrete and discrete toplogies on $X$
respectively. We use $I_{2}$ to denote the two-element indiscrete space, $%
S_{2}$ for the Sierpi\'{n}ski space (with topology $\mathcal{S}%
_{2}=\{\emptyset ,\{0\},\{0,1\}\})$ and $D_{2}$ for the two-element discrete
space. On $I_{2}$ there are only two congruences $\iota _{I_{2}}$ and $%
\upsilon _{I_{2}}$ and on $S_{2}$ there are three namely $\iota
_{S_{2}},(\Bumpeq ,\mathcal{I}_{2})$ and $\upsilon _{S_{2}}$. The space $%
D_{2}$ has five congruences $\iota _{D_{2}},(\Bumpeq ,\mathcal{S}%
_{2}),(\Bumpeq ,\{\emptyset ,\{1\},\{0,1\}\}),(\Bumpeq ,\mathcal{I}_{2})$
and $\upsilon _{D_{2}}.$ Recall, an object $Q$ in a category is called 
\textit{injective} if for any given morphism $g:C\rightarrow Q$ and
monomorphism $f:C\rightarrow B$ there exists a morphism $h:B\rightarrow Q$
such that $h\circ f=g.$ It is known (in any case easy to prove) that a
topological space is injective in the category of all topological spaces
precisely when it is an indiscrete space.

\begin{definition}
Let $\rho $ be a Hoehnke radical of topological spaces. Then $\rho $ is
called:

(1) $r$-hereditary if for every space $X$ and subspace $Y,$ $\rho _{X}\cap
Y\sqsubseteq \rho _{Y}.$

(2) $s$-hereditary if for every space $X$ and subspace $Y,$ $\rho
_{Y}\sqsubseteq \rho _{X}\cap Y.$

(3) Ideal-hereditary if it is both $r$-hereditary and $s$-hereditary.

(4) A hereditary torsion theory if it is an ideal-hereditary Hoehnke radical.
\end{definition}

\noindent Then:

\begin{proposition}
Let $\rho $ be an ideal-hereditary Hoehnke radical (= hereditary torsion
theory) in a universal class of topological spaces. Then $\rho $ is
idempotent, complete and both the radical class $\mathcal{R}_{\rho }$ and
the semisimple class $\mathcal{S}_{\rho }$ are hereditary.
\end{proposition}

\noindent For all the well-known classes of algebras, any $\rho $ as above
(ideal-hereditary Hoehnke radical) will be a Kurosh-Amitsur radical. For
topological spaces, this need not be the case. To conclude this section, we
give all the ideal-hereditary Hoehnke radicals $\rho $ of topological
spaces. There are exactly five such radicals of which three are
Kurosh-Amitsur radicals, i.e. for such $\rho $, the classes $\mathcal{R}%
_{\rho }$ and $\mathcal{S}_{\rho }$ form a corresponding pair of
connectednesses and disconnectednesses. This sounds better than it actually
is, of these three, two are trivial with the connectedness and
disconnectedness coinciding with the class of all spaces respectively. More
details about the radical-theoretic properties of the classes of topological
spaces in the next result can be found in [31].

\begin{theorem}
Let $\rho $ be an ideal-herediary Hoehnke radical of topological spaces.
Then $\rho $ is one of the following five radicals:

\noindent (a) $\rho _{X}=\upsilon _{X}$ for all $X.$ This is a KA-radical
with $\mathcal{R}_{\rho }$ the class of all spaces and $\mathcal{S}_{\rho
}=\{T\}.$

\noindent (b) $\rho _{X}=(\sim _{X},\mathbb{T}_{X})$ is the congruence with $%
x\sim _{X}y$ iff for any $U\subseteq X$ open$,$ $x\in U\Leftrightarrow y\in
U $ and $\mathbb{T}_{X}=\mathcal{T.}$ This is a KA-radical with $\mathcal{R}%
_{\rho }=\{X\mid X$ is an indiscrete space\} and $\mathcal{S}_{\rho
}=\{X\mid X$ is a $T_{0}$-space$\}.$

\noindent (c) $\rho _{X}=\iota _{X}$ for all $X.$ This is a KA-radical with $%
\mathcal{R}_{\rho }=\{T\}$ and $\mathcal{S}_{\rho }$ is the class of all
spaces.

\noindent (d) $\rho _{X}=(\Bumpeq ,\mathcal{I}_{X})$ for all $X$ which is
not a KA-radical. Here $\mathcal{R}_{\rho }=\{T\}$ and $\mathcal{S}_{\rho
}=\{X\mid X$ is an indiscrete space\}.

\noindent (e) $\rho _{X}=(\Bumpeq ,\mathbb{T}_{X})$ where $\mathbb{T}%
_{X}=\left\{ 
\begin{array}{l}
\mathcal{T}\text{ if }(X,\mathcal{T})\text{ is a }T_{1}\text{-space} \\ 
\mathcal{I}_{X}\text{ otherwise}%
\end{array}%
\right. .$ This is not KA-radical, $\mathcal{R}_{\rho }=\{T\}$ and $\mathcal{%
S}_{\rho }=\{X\mid X$ is an indiscrete space or a $T_{1}$-space\}.
\end{theorem}

\noindent In the next section we shall see that the radical-theoretic
properties presented here for topological spaces are also valid for graphs
that can have loops.

\section{Graphs that admit loops}

\noindent As for topological spaces in the previous section, here we present
the congruence approach to the connectednesses and disconnectednesses of
graphs that allow loops. Connectednesses and disconnectednesses for such
graphs have been definded and developed by Fried and Wiegandt [14] showing
that they are just the KA-radical and -semisimple classes respectively in
this category. Here it will be shown that they can be obtained as H-radicals
and then they will be characterized using congruences with conditions
resembling the classical algebraic conditions for characterizing radical
classes of rings. We have seen in the previous section that there are
hereditary torsion theories of topological spaces that can be added to those
of $S$-acts as examples of ideal-hereditary $H$-radicals which are not $KA$%
-radicals. Graphs provide a third example. In fact, there are exactly eight
ideal-hereditary H-radicals for graphs and of these, only three are
KA-radicals. Details of the work presented here can be found in [4, 33, 34].

\bigskip

\subsection{Congruences}

\noindent We start with some graph theoretic preliminaries. A graph $G$ with
vertex set $V$ and edge set $E$ will typically be denoted by $%
G=(V_{G},E_{G}),$ often without the subscripts. When we write $a\in G,$ it
actually means $a$ is a vertex of $G,$ i.e., $a\in V_{G}.$ By a graph, we
mean a non-empty vertex set, edges are not directed, no multiple edges are
allowed but loops are. For $a,b\in G,$ an edge between $a$ and $b$ is
written as $ab$ and $aa$ is the loop at $a.$ The set of all possible edges
on a graph $G$ is denoted by $C_{G}:=\{ab\mid a,b\in V_{G}\}.$ A \textit{%
(graph) homomorphism} is an edge preserving mapping from the vertex set of a
graph into the vertex set of a graph. A \textit{strong homomorphism} is a
homomorphism that sends "no edges" to "no edges" and if it is also a
bijection, it is called an \textit{isomorphism}. Isomorphic graphs $G$ and $%
H $ will be denoted by $G\cong H.$ For a graph $G=(V_{G},E_{G}),$ a \textit{%
subgraph} $H=(V_{H},E_{H})$ of $G$ is a graph with $V_{H}\subseteq V_{G}$
and $E_{H}\subseteq \{ab\mid a,b\in V_{H}$ and $ab\in E_{G}\}.$ When $%
E_{H}=\{ab\mid a,b\in V_{H}$ and $ab\in E_{G}\},$ then $H$ is called an 
\textit{induced subgraph} of $G.$ For a homomorphism $f:G\rightarrow H,$ the 
\textit{image graph} $f(G)$ will always be the induced subgraph of $H$ on
the vertex set $f(V_{G}).$ In general, unless mentioned otherwise, if a
subset $V_{H}$ of $V_{G}$ is regarded as a graph, it will be the subgraph
induced by $G$ on $V_{H}.$ There are two (non-isomorphic) one-vertex graphs,
called the\textit{\ trivial graphs}; the one with a loop $T_{0}$ and the one
without a loop $T.$ For an equivalence relation $\sim $ on a vertex set $V,$
we use $[a]$ to denote the equivalence class of $a\in V.$ For $A,B\subseteq
V_{G}$, $AB$ is the set $AB=\{ab\mid a\in A,b\in B\}\subseteq \mathcal{C}%
_{G}.$ In particular, $[a][b]=\{st\mid s\in \lbrack a],t\in \lbrack b]\}$.

\bigskip

\begin{definition}
Let $G=(V_{G},E_{G})$ be a graph. A \textit{congruence on} $G$ is a pair $%
\theta =(\sim ,\mathcal{E})$ which fulfills the following three conditions:%
\newline
(i) $\sim $ is an equivalence relation on $V_{G}$.\newline
(ii) $\mathcal{E}$ is a set of unordered pairs of elements from $V_{G},$
called the congruence edge-set, with $E_{G}\subseteq \mathcal{E\subseteq C}%
_{G}\mathcal{.}$\newline
(iii) (\textit{Substitution Property} of $\mathcal{E}$ with respect to $\sim
)$ for $x,y\in V_{G}$, if $xy\in \mathcal{E}$, then $[x][y]\subseteq 
\mathcal{E}$.\newline
A\textit{\ strong congruence on }$G$ is a pair $\theta =(\sim ,\mathcal{%
E(\sim )})$ where $\sim $ is an equivalence relation on $V_{G}$ and $%
\mathcal{E(\sim )}=\{xy\mid x,y\in V_{G}$ and $[x][y]\cap E_{G}\neq
\emptyset \}.$
\end{definition}

\medskip

\noindent Congruences are partially ordered by the relation "contained in":
for two congruences $\alpha =(\sim _{\alpha },\mathcal{E}_{\alpha })$ and $%
\beta =(\sim _{\beta },\mathcal{E}_{\beta })$ on $G,$ $\alpha $ is \textit{%
contained in }$\beta ,$ written as $\alpha \subseteq \beta ,$ if $\sim
_{\alpha }\subseteq $ $\sim _{\beta }$ and $\mathcal{E}_{\alpha }\subseteq 
\mathcal{E}_{\beta }.$ Let $\Bumpeq $ denote the identity relation on $V_{G}$
$($also called the diagonal; i.e., $x\Bumpeq y$ if and only if $x=y).$ The
congruence $\iota _{G}:=(\Bumpeq ,E_{G})$ is called the \textit{identity
congruence} on $G$ and is the smallest congruence on $G.$ The\textit{\
universal congruence} on $G$ is the pair $\upsilon
_{G}=(\leftrightsquigarrow ,\mathcal{E})$ where $\leftrightsquigarrow $ is
the universal relation $($i.e., $a\leftrightsquigarrow b$ for all $a,b\in
V_{G})$ and $\mathcal{E}$ $=\{ab\mid a,b\in V_{G}\}.$ Any congruence on $G$
contains $\iota _{G}$ and is contained in $\upsilon _{G}.$ Given any graph
homomorphism $f:G\longrightarrow H$, the \textit{kernel of} $f$, written as $%
\ker f=(\sim _{f},\mathcal{E}_{f}),$ is the congruence with $\sim
_{f}=\{(x,y)\mid x,y\in V_{G},f(x)=f(y)\}$ and $\mathcal{E}_{f}=\{uv\mid
u,v\in V_{G},f(u)f(v)\in E_{H}\}.$ With $f$ is also associated the \textit{%
strong kernel of }$f,$ defined by \textit{sker} $f=(\sim _{f},\mathcal{%
E(\sim }_{f}\mathcal{)}).$ This is a strong congruence on $G$ and \textit{%
sker} $f\subseteq \ker f;$ in fact, if $\theta =(\sim _{f},\mathcal{E})$ is
any congruence on $G$ for some $\mathcal{E}$, then \textit{sker} $f\subseteq
\theta .$ If $f$ is a strong homomorphism, then $\ker f=$ \textit{sker} $f.$
Given any congruence $\theta =(\sim ,\mathcal{E})$ on a graph $%
G=(V_{G},E_{G})$, we define a new graph, denoted by $G/\theta =(V_{G/\theta
},E_{G/\theta })$ and called \textit{the quotient of $G$ modulo} $\theta ,$
by taking $V_{G/\theta }:=\{[x]\mid x\in V_{G}\}$ and $E_{G/\theta
}:=\{[x][y]\mid xy\in \mathcal{E}\}.$ The natural or canonical mapping $%
p_{\theta }:G\rightarrow G/\theta $ given by $p_{\theta }(x)=[x]$ is a
surjective homomorphism with $\ker p_{\theta }=\theta .$ For $\theta =\iota
_{G}$ we have $G/\iota _{G}$ isomorphic to $G.$ If $\theta $ is a strong
congruence, then $p_{\theta }$ is a strong homomorphism with \textit{sker} $%
p_{\theta }=\theta .$ In general, if we take $\theta =(\Bumpeq ,\mathcal{E})$
for some suitable $\mathcal{E}$ to make $\theta $ a congruence on $G,$ then $%
G/\theta $ is the graph with vertex set $V_{G/\theta }=V_{G}$ and edge set $%
E_{G/\theta }=\mathcal{E}$ (here we identify $[x]=\{x\}$ with $x).$ If $%
\upsilon _{G}$ is the universal congruence on $G,$ then $G/\upsilon _{G}$ is
isomorphic to the trivial graph $T_{0}$ with a loop.

\bigskip

\noindent \textbf{Ordering of congruences. }For a given graph $G,$ we denote
the set of all congruences on $G$ by $Con(G)$ which is a partially ordered
set with respect to containment $\subseteq .$ But we can say more. $Con(G)$
is a bounded complete lattice with the meet and join of $\theta _{i}\in
Con(G),i\in I,$ given by the intersection $\bigcap\limits_{i\in I}\theta
_{i} $ and sum $\sum\limits_{i\in I}\theta _{i}$ respectively. These are
defined as follows: for any collection of congruences $\{\theta _{i}=(\sim
_{i},\mathcal{E}_{i})\mid i\in I\}\subseteq Con(G)$ the greatest lower bound
in $Con(G)$ is given by $\bigcap\limits_{i\in I}\theta _{i}=(\sim _{\cap },%
\mathcal{E}_{\cap })$ where $a\sim _{\cap }b\Leftrightarrow a\sim _{i}b$ for
all $i\in I$ and $ab\in \mathcal{E}_{\cap }\Leftrightarrow ab\in \mathcal{E}%
_{i}$ for all $i\in I.$ The smallest upper bound is $\sum\limits_{i\in
I}\theta _{i}=(\sim _{\Sigma },\mathcal{E}_{\Sigma })$ where $a\sim _{\Sigma
}b\Leftrightarrow \exists $ $c_{1},c_{2},...,c_{n}\in G$ with $a=c_{1}\sim
_{i_{1}}c_{2}\sim _{i_{2}}c_{3}\sim _{i_{3}}...\sim _{i_{n-1}}c_{n}=b$ for
some $i_{1},i_{2},...,i_{n-1}\in I,n\geq 2$ and $\mathcal{E}_{\Sigma
}=\{ab\mid ab\in \lbrack x]_{\Sigma }[y]_{\Sigma }$ for some $xy\in
\bigcup\limits_{i\in I}\mathcal{E}_{i}\}.$ If $\theta _{i}$ is strong for
all $i\in I,$ then so is $\sum\limits_{i\in I}\theta _{i}.$ For two
congruences $\alpha $ and $\beta ,$ we write the sum as $\alpha +\beta .$

\bigskip

\noindent \textbf{Image of a congruence. }Let $f:G\rightarrow H$ be a
surjective homomorphism with $\theta =(\sim _{\theta },\mathcal{E}_{\theta
}) $ a congruence on $G$ and $\alpha =\ker f.$ The\textit{\ image of }$%
\theta $\textit{\ under }$f$ is the congruence $f(\theta )=(\sim _{f(\theta
)},\mathcal{E}_{f(\theta )})$ of $H$ where $f(a)\sim _{f(\theta
)}f(b)\Leftrightarrow a\sim _{\alpha +\theta }b.$ Recall that $x\sim _{\ker
f}y\Leftrightarrow f(x)=f(y);$ hence $f(a)\sim _{f(\theta
)}f(b)\Leftrightarrow \exists $ $c_{1},c_{2},...,c_{n}\in G$ with $%
f(a)=f(c_{1}),c_{1}\sim _{\theta }c_{2},f(c_{2})=f(c_{3}),c_{3}\sim _{\theta
}c_{4},...,f(c_{n-1})=f(c_{n}),c_{n}\sim _{\theta }b.$ Let $\mathcal{E}%
_{f(\theta )}=\{cd\mid cd\in \lbrack f(x)]_{f(\theta )}[f(y)]_{f(\theta )}$
for some $x,y\in G$ with $xy\in \mathcal{E}_{\theta }$ or $f(x)f(y)\in
E_{H}\}.$ It can be checked that $\sim _{f(\theta )}$ is well-defined and in
the chain of equalities and equivalences, it does not matter whether one
starts or ends with an equality or equivalence. Then $f(\theta )$ is a
congruence which is strong if $\theta $ is strong.

\bigskip

\noindent \textbf{Isomorphism theorems for congruences.}

\begin{theorem}
(\textit{First Isomorphism Theorem}) Let $f:G\rightarrow H$ be a surjective
homomorphism. Then $G/\ker f$ is isomorphic to $H$.
\end{theorem}

\noindent Let $G$ be a graph with induced subgraph $H.$ Then a congruence $%
\theta =(\sim ,\mathcal{E})$ on $G$ induces a congruence $H\cap \theta
=(\sim _{H},\mathcal{E}_{H})$ on $H$ with $\sim _{H}=(V_{H}\times V_{H})\cap 
$ $\sim $ $=\{(a,b)\mid a,b\in V_{H}$ and $a\sim b\}$ and $\mathcal{E}%
_{H}=\{ab\mid a,b\in V_{H}\}\cap \mathcal{E}=\{ab\mid a,b\in V_{H}$ with $%
ab\in \mathcal{E}\}.$ The mapping $f:H\rightarrow G/\theta $ defined by $%
f(a)=[a]$ for all $a\in V_{H}$ is a homomorphism with $\ker f=H\cap \theta .$
Now $f(V_{H})$ is a set of vertices of $G/\theta $ on which we form the
induced subgraph of $G/\theta $, denoted by $(H+\theta )/\theta .$ Then, by
the First Isomorphism Theorem, we have: \medskip

\begin{theorem}
(\textit{Second Isomorphism Theorem}) Let $H$ be an induced subgraph of a
graph $G.$ Let $\theta $ be a congruence on $G.$ Then $H\cap \theta $ as
defined above is a congruence on $H$ and $H/H\cap \theta \cong (H+\theta
)/\theta $ where $(H+\theta )/\theta $ is the induced subgraph of $G/\theta $
on the vertex set $\{[a]\mid a\in V_{H}\}.$
\end{theorem}

\begin{theorem}
(\textit{Third Isomorphism Theorem}) Let $G$ be a graph with two congruences 
$\theta _{1}=(\sim _{1},\mathcal{E}_{1})$ and $\theta _{2}=(\sim _{2},%
\mathcal{E}_{2})$ on $G$ for which $\theta _{1}\subseteq \theta _{2}.$ Then $%
\theta _{2}/\theta _{1}:=(\sim ,\mathcal{E})$ is a congruence on $G/\theta
_{1}$ where $[a]_{1}\sim \lbrack b]_{1}\Leftrightarrow a\sim _{2}b$ and $%
[a]_{1}[b]_{1}\in \mathcal{E}\Leftrightarrow ab\in \mathcal{E}_{2}.$
Moreover, $(G/\theta _{1})/(\theta _{2}/\theta _{1})$ is isomorphic to $%
G/\theta _{2}.$
\end{theorem}

\begin{corollary}
Let $G$ be a graph with $\theta $ a fixed congruence on $G.$ Any congruence $%
\xi $ on the graph $G/\theta $ is of the form $\alpha /\theta $ for some
congruence $\alpha $ on $G$ with $\theta \subseteq \alpha .$ Moreover, there
is a one-to-one correspondence between $\{\alpha \mid \alpha $ is a
congruence on $G$ with $\theta \subseteq \alpha \}$ and $Con(G/\theta )$
which preserves inclusions, intersections and unions of congruences.
\end{corollary}

\bigskip

\noindent \textbf{Subdirect product of graphs. }For an index set $I,$ let $%
G_{i}=(V_{i},E_{i})$ be a graph for all $i\in I.$ The product $%
\prod\limits_{i\in I}G_{i}$ of the graphs $G_{i}$ is the graph $%
\prod\limits_{i\in I}G_{i}:=(\prod\limits_{i\in I}V_{i},E)$ where $%
\prod\limits_{i\in I}V_{i}$ is just the usual Cartesian product of the sets $%
V_{i}$ and $E=\{fg\mid f,g\in \prod\limits_{i\in I}V_{i}$ with $f(i)g(i)\in
E_{i}$ for all $i\in I\}.$ For every $j\in I,$ the $j$-th \textit{projection}
$\pi _{j}:\prod\limits_{i\in I}G_{i}\rightarrow G_{j}$ defined by $\pi
_{j}(f)=f(j)$ for all $f\in \prod\limits_{i\in I}V_{i}$ is a surjective
homomorphism. An induced subgraph $H$ of $\prod\limits_{i\in I}G_{i}$ is
called a \textit{subdirect product of the graphs} $G_{i},i\in I,$ provided
the restriction of each projection $\pi _{j}$ to $H$ is a surjective mapping
onto $G_{j}$. As in universal algebra, subdirect products can be
characterized in terms of congruences and quotients:

\begin{theorem}
A graph $G$ is a subdirect product of graphs $G_{i},i\in I,$ if and only if
for every $i\in I$ there are congruences $\theta _{i}$ on $G$ with $G_{i}$
isomorphic to $G/\theta _{i}$ and $\bigcap\limits_{i\in I}\theta _{i}=\iota
_{G}.$
\end{theorem}

\bigskip

\subsection{Radical theory}

\noindent All radical theoretic considerations are in a universal class of
graphs $\mathcal{W.}$ This means $\mathcal{W}$ is non-empty, closed under
homomorphic images and closed under the taking of subgraphs (= strongly
hereditary). We do not distinguish between isomorphic graphs. From the
definition, it follows that $\mathcal{W}$ contains a one-vertex graph, and
consequently also the class $\mathcal{T}=\{T_{0},T\}$ of all trivial graphs
is contained in $\mathcal{W}$. Note that there is a unique congruence on $%
T_{0}$ since here $\iota _{T_{0}}=\upsilon _{T_{0}},$ but on $T$ these two
congruences are distinct. A class of graphs $\mathcal{M}$ in $\mathcal{W}$
is an \textit{abstract class} provided it is closed under isomorphic copies
and it contains the trivial graph $T_{0}$. All subclasses of $\mathcal{W}$
under consideration will be assumed to be abstract, even though it may not
always be explicitly stated. Since $G/\upsilon _{G}\cong T_{0}$ for any
graph $G,$ there is always at least one congruence $\theta $ on a graph $G$
for which $G/\theta $ is in any abstract class. For a class $\mathcal{M}$ in 
$\mathcal{W},$ we use $\overline{\mathcal{M}}$ to denote the \textit{%
subdirect closure} of $\mathcal{M},$ i.e., the class of all graphs that are
subdirect products of graphs from $\mathcal{M}.$ Clearly $\mathcal{M}%
\subseteq \overline{\mathcal{M}}$ and we say $\mathcal{M}$ is \textit{%
subdirectly closed} if $\mathcal{M}=\overline{\mathcal{M}}.$ We start with
the definition of the Hoehnke radical for which we will need the image of a
congruence under a homomorphism as defined in the previous section.

\begin{definition}
An\textit{\ H-radical} on $\mathcal{W}$ is a function $\varrho $ that
assigns to every graph $G$ in $\mathcal{W}$ a congruence $\varrho
(G)=\varrho _{G}$ on $G$ such that:

(H1) For every surjective homomorphism $f:G\rightarrow H,$ $f(\varrho
(G))\subseteq \varrho (H).$

(H2) For any graph $G$, $\varrho (G/\varrho _{G})=\iota _{G/\varrho _{G}}$,
the identity congruence on $G/\varrho _{G}$.
\end{definition}

\noindent For an H-radical $\varrho ,$ if $\varrho (G)=\iota _{G},$ then $G$
is called \textit{semisimple} (actually, $\varrho $-semisimple), the class $%
\mathcal{S}_{\varrho }=\{G\in \mathcal{W}\mid \varrho (G)=\iota _{G}\}$ is
called the \textit{associated semisimple class} and $\mathcal{R}_{\varrho
}=\{G\in \mathcal{W}\mid G/\varrho _{G}$ is a trivial graph\} is the \textit{%
associated} \textit{radical class}. For an H-radical $\varrho $ on $\mathcal{%
W},$ note that $\{T_{0}\}\subseteq \mathcal{S}_{\varrho }\cap \mathcal{R}%
_{\varrho }\subseteq \mathcal{T}\subseteq $ $\mathcal{R}_{\varrho }$ and the
radical class $\mathcal{R}_{\varrho }$ is always homomorphically closed. If $%
G/\varrho _{G}$ is trivial, then $\varrho _{G}=(\leftrightsquigarrow ,%
\mathcal{C}_{G})$ for $E_{G}\neq \emptyset $ and when $E_{G}=\emptyset ,$
then $\varrho _{G}$ can be $(\leftrightsquigarrow ,\mathcal{C}_{G})$ or $%
(\leftrightsquigarrow ,\emptyset ).$ An H-radical is very general and is
always of a prescribed form as the next result shows.

\begin{theorem}
Let $\varrho $ be a mapping that assigns to any graph $G$ in $\mathcal{W}$ a
congruence $\varrho (G)=\varrho _{G}$ on $G.$ Then $\varrho $ is a H-radical
on $\mathcal{W}$ if and only if there is an abstract class of graphs $%
\mathcal{M}$ in $\mathcal{W}$ such that for all $G$ in $\mathcal{W},$ $%
\varrho (G)=\cap \{\theta \mid \theta $ is a congruence on $G$ for which $%
G/\theta \in \mathcal{M}\}.$ Furthermore, $\mathcal{S}_{\varrho }=$ $%
\overline{\mathcal{M}}.$
\end{theorem}

\noindent The salient properties of an H-radical are contained in the next
corollary.

\begin{corollary}
(1) The semisimple class of any H-radical is subdirectly closed. \newline
(2) For an H-radical $\varrho ,\varrho (G)$ is the smallest congruence on $G$
for which $G/\varrho (G)$ is semisimple $($i.e., if $\theta $ is a
congruence on $G$ with $G/\theta \in \mathcal{S}_{\varrho },$ then $\varrho
(G)\subseteq \theta ).$ Or, equivalently, $G/\varrho (G)$ is the largest
semisimple image of $G$ (in the following sense: if $g:G\rightarrow H$ is a
surjective homomorphism with $H\in \mathcal{S}_{\varrho },$ then there is a
homomorphism $h:G/\varrho (G)\rightarrow H$ such that $h\circ p=g$ where $%
p:G\rightarrow G/\varrho (G)$ is the canonical quotient map). \newline
(3) For any abstract class of graphs $\mathcal{M}$ in $\mathcal{W}$ and $%
G\in \mathcal{W}$, define $\varrho (G):=\cap \{\theta \mid \theta \in Con(G)$
with $G/\theta \in \mathcal{M}\}.$ Then $\varrho $ is an H-radical with $%
\mathcal{S}_{\varrho }=$ $\overline{\mathcal{M}},$ i.e., every semisimple
graph is a subdirect product of graphs from $\mathcal{M}.$
\end{corollary}

\noindent A class $\mathcal{M}$ of graphs is said to be \textit{hereditary}
(respectively \textit{strongly hereditary}) if $G\in \mathcal{M}$ implies
all the induced subgraphs of $G$ (respectively all the subgraphs of $G)$ are
in $\mathcal{M}$. Hereditariness is retained under subdirect closure. Next
we recall the definitions of Kurosh-Amitsur radical and semisimple classes
of graphs. As for topological spaces, these classes of graphs are called
connectednesses and disconnectednesses respectively. The terminology used is
obvious when looking at the many examples of these classes ([14]). A class $%
\mathcal{C}\subseteq \mathcal{W}$ is a \textit{connectedness} (= \textit{%
KA-radical class)} if it satisfies the following condition: A graph $G\in 
\mathcal{W}$ is in $\mathcal{C}$ if and only if every non-trivial
homomorphic image of $G$ has a non-trivial induced subgraph which is in $%
\mathcal{C}.$ A class $\mathcal{D}\subseteq \mathcal{W}$ is a \textit{%
disconnectedness (= KA-semisimple class)} if it satisfies the following
condition: A graph $G\in \mathcal{W}$ is in $\mathcal{D}$ if and only if
every non-trivial induced subgraph of $G$ has a non-trivial homomorphic
image which is in $\mathcal{D}.$ For examples and many additional statements
and properties of these classes of graphs, see Fried and Wiegandt [14]. The
class of trivial graphs $\mathcal{T}$ is always contained in any
connectedness and also in any disconnectedness. It is easy to find examples
of connectednesses and disconnectednesses: If $\mathcal{M}\subseteq \mathcal{%
W}$ is a hereditary class, then $\mathcal{UM}:=\{G\in \mathcal{W\mid G}$ has
no non-trivial homomorphic image in $\mathcal{M}\}$ is a connectedness and
if $\mathcal{H}\subseteq \mathcal{W}$ is a homomorphically closed class,
then $\mathcal{SH}:=\{G\in \mathcal{W}\mid G$ has no non-trivial induced
subgraph in $\mathcal{H}\}$ is a disconnectedness. From the preceding, we
thus have: If $\mathcal{C}$ is a connectedness, then $\mathcal{SC}$ is a
disconnectedness and if $\mathcal{D}$ is a disconnectedness, then $\mathcal{%
UD}$ is a connectedness. Moreover, it can be shown that a class $\mathcal{C}%
\subseteq \mathcal{W}$ is a connectedness if and only if $\mathcal{C}=%
\mathcal{USC}$ and a class $\mathcal{D}\subseteq \mathcal{W}$ is a
disconnectedness if and only if $\mathcal{D}=\mathcal{SUD}.$ If $\varrho $
is a Hoehnke radical, then $\mathcal{R}_{\varrho }=\mathcal{US}_{\varrho }$
and if $\mathcal{S}_{\varrho }$ is hereditary, then $\mathcal{S}_{\varrho
}\subseteq \mathcal{SR}_{\varrho }.$ Two properties that a Hoehnke radical $%
\varrho $ on $\mathcal{W}$ may satisfy are:

\noindent \textit{Complete}: If $\theta $ is a strong congruence on $G\in 
\mathcal{W}$ with $[a]_{\theta }\in \mathcal{R}_{\varrho }$ for all $a\in
V_{G},$ then $\theta \subseteq \varrho _{G};$ and

\noindent \textit{Idempotent}: For $G\in \mathcal{W}$ and all $a\in V_{G},$ $%
[a]_{\varrho _{G}}\in \mathcal{R}_{\varrho }.$ \newline
Then we have:

\begin{theorem}
Let $\varrho $ be an H-radical on $\mathcal{W}$ which is complete,
idempotent and such that for all $G\in \mathcal{W},$ $\varrho _{G}$ is a
strong congruence on $G.$ Then $\mathcal{S}_{\varrho }$ is a
disconnectedness and $\mathcal{R}_{\varrho }=\mathcal{US}_{\varrho }$ is a
connectedness. Conversely, suppose $\mathcal{D}$ is a disconnectedness in $%
\mathcal{W}$ with corresponding connectedness $\mathcal{C}.$ Then there is
an H-radical $\varrho $ on $\mathcal{W}$ which is complete, idempotent and
for all $G\in \mathcal{W},$ $\varrho _{G}$ is a strong congruence on $G.$
Moreover, $\mathcal{S}_{\varrho }=\mathcal{D}$ and $\mathcal{R}_{\varrho }=%
\mathcal{UD}=\mathcal{C.}$
\end{theorem}

\noindent This result motivates the following terminology: An H-radical $%
\varrho $ which is complete, idempotent and for which $\varrho _{G}$ is a
strong congruence for all $G$ is called a \textit{KA-radical}. Next we show
that a connectedness of graphs can be characterized with conditions using
congruences which correspond to the conditions characterizing a radical
class of associative rings using ideals. In particular, it is shown that any
graph has a largest radical congruence which contains all other radical
congruences on the graph and the quotient of a graph by this largest radical
congruence has no non-trivial radical congruences. A significant property of
KA-semisimple classes in general radical theory is that any object in the
universal class has a maximal semisimple image. This is also the case for
disconnectednesses of graphs as Fried and Wiegandt have shown [14]. This
result is recalled below, together with their characterization of
connectednesses. Remember, whenever a subset of vertices of a graph $G$ is
considered as a graph and nothing else is mentioned, it is the subgraph
induced by the graph $G.$

\begin{theorem}
\lbrack 14] Let $\mathcal{C}$ be a connectedness with corresponding
disconnectedness $\mathcal{D}=\mathcal{SC}.$ Then: \newline
(1) For every $G\in \mathcal{W,}$ there is a strong homomorphism $%
q:G\rightarrow G_{\mathcal{D}}$ with $G_{\mathcal{D}}\in \mathcal{D}$ and if 
$f:G\rightarrow H$ is any surjective homomorphism with $H\in \mathcal{D},$
then there is a homomorphism $g:G_{\mathcal{D}}\rightarrow H$ such that $%
g\circ q=f$. (In categorical terms, this means $\mathcal{D}$ is an
epi-reflective subcategory of $\mathcal{W}.$) $G_{\mathcal{D}}$ is called
the \textit{maximal }$D$\textit{-image }of $G.$ \newline
(2) For every $a\in G_{\mathcal{D}},$ $q^{-1}(a)\in \mathcal{C}$ and it is
maximal in the sense that it is not properly contained in any other induced
subgraph of $G$ which is in $\mathcal{C}.$ \newline
(3) If $H$ is an induced subgraph of $G$ with $H\in \mathcal{C},$ then there
is an $a\in G_{\mathcal{D}}$ such that $H\subseteq q^{-1}(a).$
\end{theorem}

\noindent Let $\mathcal{C}$ be an abstract subclass of $\mathcal{W.}$ $%
\mathcal{C}$ \textit{closed under forming strings} means if a graph $G$ is
the union of induced subgraphs $G_{i},i\in I$ with each $G_{i}\in \mathcal{C}
$ and for every $j,k\in I$ there are finitely many indices $%
j=i_{1},i_{2},...,i_{n}=k$ in $I$ with $G_{i_{t-1}}\cap G_{i_{t}}\neq
\emptyset $ for $t=2,3,...,n,$ then $G\in \mathcal{C};$ and $\mathcal{C}$ 
\textit{weakly extensive} means whenever $f:G\rightarrow H$ is a strong
homomorphism with $H$ and $f^{-1}(b)$ in $\mathcal{C}$ for all $b\in H,$
then $G\in \mathcal{C}.$ The intersection $G_{i_{t-1}}\cap G_{i_{t}}$ refers
to the vertex sets of the graphs.

\begin{theorem}
\lbrack 14] Let $\mathcal{C}$ be an abstract subclass of $\mathcal{W}.$ Then 
$\mathcal{C}$ is a connectedness if and only if it satisfies the following
three conditions:

\noindent (1) $\mathcal{C}$ is homomorphically closed.

\noindent (2) $\mathcal{C}$ is closed under forming strings.

\noindent (3) $\mathcal{C}$ is weakly extensive.
\end{theorem}

\noindent For an abstract subclass $\mathcal{C}$ of $\mathcal{W}$, \textit{a}
subgraph $H$ of $G$ is a $\mathcal{C}$\textit{-subgraph} of $G$ provided $%
H\in \mathcal{C;}$ likewise an \textit{induced }$\mathcal{C}$\textit{%
-subgraph }is defined. A congruence $\alpha $ on a graph $G$ is a $\mathcal{C%
}$\textit{-congruence on }$G$ if $\alpha $ is a strong congruence and $[a]$
is an induced $\mathcal{C}$-subgraph of $G$ for all $a\in G.$ A congruence
on a graph $G$ is called \textit{trivial} if $[a]=\{a\}$ for all $a\in G.$
Note that for any $\mathcal{E}$ with $E_{G}\subseteq \mathcal{E}\subseteq
C_{G}$, $\alpha =(\Bumpeq ,\mathcal{E})$ is a trivial congruence on $G,$ but
it need not coincide with the identity congruence $\iota _{G}:=(\Bumpeq
,E_{G})$ on $G.$ If however, $\alpha $ is a trivial strong congruence, then $%
\alpha =\iota _{G}$ since in this case, $\mathcal{E}(\Bumpeq )=E_{G}$. For
our characterization of connectednesses, we start with:

\begin{proposition}
Let $\mathcal{C}$ be an abstract subclass of $\mathcal{W}.$ Then $\mathcal{C}
$ is a connectedness if and only if it satisfies the following condition: A
graph $G\in \mathcal{W}$ is in $\mathcal{C}$ $\Leftrightarrow $ every
non-trivial homomorphic image of $G$ has a non-trivial $\mathcal{C}$%
-congruence.
\end{proposition}

\noindent Let $\mathcal{C}$ be an abstract class of graphs in $\mathcal{W.}$
For any graph $G,$ define a congruence $\varrho (G)$ on $G$ by $\varrho
(G)=\sum \{\alpha \mid \alpha $ is a $\mathcal{C}$-congruence on $G\}.$ This
sum is not void, since $\iota _{G}$ is always a $\mathcal{C}$-congruence on $%
G.$ Moreover, $\varrho (G)$ is a strong congruence on $G$ and it contains
all $\mathcal{C}$-congruences on $G.$ Usually we write $\varrho (G)$ as $%
\varrho _{G},$ but there are instances when the former will be better to
use. Note that $\upsilon _{G}$ is a $\mathcal{C}$-congruence on $%
G\Leftrightarrow G\in \mathcal{C}.$

\begin{proposition}
Let $\mathcal{C}$ be a connectedness in $\mathcal{W}$ with associated
disconnectedness $\mathcal{D}=\mathcal{SC}.$ For any graph $G$ in $W,$ $%
\varrho _{G}=\sum \{\alpha \mid \alpha $ is a $\mathcal{C}$-congruence on $%
G\}$ is a $\mathcal{C}$-congruence on $G$ which contains all $\mathcal{C}$%
-congruences on $G$ and $\varrho _{G}=\sigma _{G}$ where $\sigma _{G}=\cap
\{\theta \mid \theta \in Con(G)$ with $G/\theta \in \mathcal{D}\}.$
\end{proposition}

\noindent In terms of congruences and the terminology introduced in this
section, the condition weakly extensive can be rephrased as: whenever $%
\alpha $ is a $\mathcal{C}$-congruence of $G$ and $G/\alpha \in \mathcal{C},$
then $G\in \mathcal{C}.$ Furthermore, this condition can be formulated using
only congruences: The class $\mathcal{C}$ is weakly extensive if and only if
it satisfies: whenever $\alpha $ and $\beta $ are congruences on $G$ with $%
\alpha \subseteq \beta ,$ $\alpha $ a $\mathcal{C}$-congruence on $G$ and $%
\beta /\alpha $\ a $\mathcal{C}$-congruence on $G/\alpha ,$ then $\beta $ is
a $\mathcal{C}$-congruence on $G.$ One last definition: $\mathcal{C}$ is 
\textit{inductive} if it satisfies: whenever $\alpha _{1}\subseteq \alpha
_{2}\subseteq \alpha _{3}\subseteq ...$ is a chain of $\mathcal{C}$%
-congruences on the graph $G,$ then $\sum\limits_{i=1}^{\infty }\alpha _{i}$
is a $\mathcal{C}$-congruence on $G.$ We now give the main result which
shows that a connectedness of graphs can be characterized using congruences
in exactly the same way that a radical class of associative rings can be
characterized using ideals (again, compare this with Theorem 2.1).

\begin{theorem}
Let $\mathcal{C}$ be an abstract class of graphs in $\mathcal{W.}$ Then
statements (1), (2) and (3) below are equivalent:

\noindent (1) $\mathcal{C}$ is a connectedness.

\noindent (2) If for all $G\in W,$ $\varrho _{G}:=\sum \{\alpha \mid \alpha $
is a $\mathcal{C}$-congruence on $G\},$ then $\mathcal{C}$ satisfies the
following three conditions:

(a) For every surjective homomorphism $f:G\rightarrow H,$ $f(\varrho
_{G})\subseteq \varrho _{H}.$

(b) For any graph $G$ in $W,$ $\varrho _{G}$ is a $\mathcal{C}$-congruence
on $G$ which contains all $\mathcal{C}$-congruences on $G.$

(c) $\varrho (G/\varrho _{G})=\iota _{G/\varrho _{G}}$ for all graphs $G.$

\noindent (3) $\mathcal{C}$ satisfies the following three conditions:

(a) If $\alpha $ is a congruence on $G$ and $G\in \mathcal{C,}$ then $%
G/\alpha \in \mathcal{C.}$

(b) $\mathcal{C}$ is inductive.

(c) $\mathcal{C}$ is weakly extensive.
\end{theorem}

\subsection{Hereditary torsion theories\textbf{.}}

\noindent In this section, we take the universal class $\mathcal{W}$ to be
the class of all graphs.

\begin{definition}
Let $\varrho $ be a Hoehnke radical of graphs. Then $\varrho $ is called:

(1) $r$-hereditary if for every graph $G$ and induced subgraph $H$ of $G,$ $%
\varrho _{G}\cap H\subseteq \varrho _{H}.$

(2) $s$-hereditary if for every graph $G$ and induced subgraph $H$ of $G,$ $%
\varrho _{H}\subseteq \varrho _{G}\cap H.$

(3) Ideal-hereditary if it is both $r$-hereditary and $s$-hereditary.

(4) A hereditary torsion theory if it is an ideal-hereditary Hoehnke radical.
\end{definition}

\noindent We will need:

\begin{proposition}
Let $\varrho $ be a Hoehnke radical.

\noindent $(1)$ If $\varrho $ is $r$-hereditary, then $\varrho $ is
idempotent and $\mathcal{R}_{\varrho }$ is hereditary.

\noindent $(2)$ $\varrho $ is $s$-hereditary if and only if $\mathcal{S}%
_{\varrho }$ is hereditary. Any one of these two conditions, implies that $%
\varrho $ is complete.
\end{proposition}

\begin{corollary}
Let $\varrho $ be an ideal-hereditary Hoehnke radical (= hereditary torsion
theory). Then $\varrho $ is idempotent, complete and both the associated
radical class $\mathcal{R}_{\varrho }$ and the associated semisimple class $%
\mathcal{S}_{\varrho }$ are hereditary.
\end{corollary}

\noindent For all the well-known classes of algebras, any $\varrho $ as
above (ideal-hereditary H-radical) will be a KA-radical meaning the
associated radical class and semisimple class is a KA-radical class and a
KA-semisimple class respectively. For graphs, as is the case for topological
spaces (previous section) and $S$-acts (see Wiegandt [37]), this need not be
the case. Below it is shown that there are exactly eight ideal-hereditary
H-radicals of which three are KA-radicals. As will be seen, they are
determined by the six non-isomorphic two-vertex graphs which will be denoted
by $B_{i},i=1,2,3,...,6$ where $V_{B_{i}}=\{0,1\}$ and $E_{B_{1}}=\emptyset $
(empty set), $E_{B_{2}}=\{01\},E_{B_{3}}=\{00\},E_{B_{4}}=\{00,11%
\},E_{B_{5}}=\{01,11\}$ and $E_{B_{6}}=\{00,01,11\}.$ Let $\mathcal{B}$
denote the set consisting of the six two-vertex graphs and for a graph $G,$ $%
L_{G}$ will be the set of all vertices from $G$ which has a loop, i.e., $%
L_{G}=\{t\in G\mid tt\in E_{G}\}.$ Recall, for an equivalence $\sim $ on the
vertex set $V_{G},$ $(\sim ,\mathcal{E}(\sim ))$ denotes the strong
congruence on $G$ determined by $\sim .$ As was the case for topological
spaces, of the three ideal-hereditary $H$-radicals which are KA-radical,
only one is non-trivial.

\begin{theorem}
Let $\varrho $ be an ideal-hereditary Hoehnke radical. Then $\varrho $ is
one of the following eight radicals:

\noindent (a) $\varrho _{G}=(\leftrightsquigarrow ,\mathcal{E}%
(\leftrightsquigarrow ))$ for all $G.$ This is a Kurosh-Amitsur radical with 
$\mathcal{R}_{\varrho }$ the class of all graphs, $\mathcal{S}_{\varrho }=%
\mathcal{T}$ and $\mathcal{S}_{\varrho }\cap \mathcal{B=\emptyset .}$

\noindent (b) $\varrho _{G}=\nu _{G}$ for all $G.$ This is not a
Kurosh-Amitsur radical, $\mathcal{R}_{\varrho }$ is the class of all graphs, 
$\mathcal{S}_{\varrho }=\{T_{0}\}$ and $\mathcal{S}_{\varrho }\cap \mathcal{%
B=\emptyset .}$

\noindent (c) For all $G,\varrho _{G}=(\sim _{G},\mathcal{E}_{G})$ where $%
\sim _{G}$ is the equivalence relation with equivalence classes $\{\{a\}\mid
a\in V_{G}-L_{G}\}\cup \{L_{G}\}$ and $\mathcal{E}_{G}=\mathcal{E}(\sim
_{G}).$ This is a Kurosh-Amitsur radical with $\mathcal{R}_{\varrho
}=\{G\mid $ if $G$ is non-trivial, then every vertex of $G$ has a loop\}, $%
\mathcal{S}_{\varrho }=\{G\mid G$ has at most one loop\} and $\mathcal{S}%
_{\varrho }\cap \mathcal{B=\{}B_{1},B_{2},B_{3},B_{5}\}.$

\noindent (d) For all $G,\varrho _{G}=(\sim _{G},\mathcal{E}_{G})$ where $%
\sim _{G}$ is the equivalence relation with equivalence classes $\{\{a\}\mid
a\in L_{G}\}\cup \{V_{G}-L_{G}\}$ and $\mathcal{E}_{G}=E_{G}\cup \{tt\mid
t\in V_{G}-L_{G}\}.$ This is not a Kurosh-Amitsur radical, $\mathcal{R}%
_{\varrho }=\mathcal{T,}$ $\mathcal{S}_{\varrho }=\{G\mid $every vertex of $%
G $ has a loop\} and $\mathcal{S}_{\varrho }\cap \mathcal{B=\{}%
B_{4},B_{6}\}. $

\noindent (e) $\varrho _{G}=(\Bumpeq ,\mathcal{C}_{G})$ for all $G.$ This is
not a Kurosh-Amitsur radical, $\mathcal{R}_{\varrho }=\mathcal{T,}$ $%
\mathcal{S}_{\varrho }$ is the class of all complete graphs with a loop at
every vertex and $\mathcal{S}_{\varrho }\cap \mathcal{B=\{}B_{6}\}.$

\noindent (f) $\varrho _{G}=\iota _{G}$ for all $G.$ This is a
Kurosh-Amitsur radical with $\mathcal{R}_{\varrho }=\mathcal{T,}$ $\mathcal{S%
}_{\varrho }$ the class of all graphs and $\mathcal{B\subseteq \mathcal{S}%
_{\varrho }.}$

\noindent (g) For all $G,\varrho _{G}=(\Bumpeq ,\mathcal{E}_{G})$ where $%
\mathcal{E}_{G}=E_{G}\cup \{ab\mid a,b\in L_{G}\}.$ This is not a
Kurosh-Amitsur radical, $\mathcal{R}_{\varrho }=\mathcal{T,}$ $\mathcal{S}%
_{\varrho }=\{G\mid $ if $G$ is non-trivial, then $ab\in E_{G}$ for all $%
a,b\in L_{G}\}$ and $\mathcal{S}_{\varrho }\cap \mathcal{B=\{}%
B_{1},B_{2},B_{3},B_{5},B_{6}\}.$

\noindent (h) For all $G,\varrho _{G}=(\Bumpeq ,\mathcal{E}_{G})$ where $%
\mathcal{E}_{G}=E_{G}\cup L_{G}V_{G},\mathcal{R}_{\varrho }=\mathcal{T,}$ $%
\mathcal{S}_{\varrho }=\{G\mid $ if $a\in L_{G},$ then $ab\in E_{G}$ for all 
$b\in V_{G}\}$ and $\mathcal{S}_{\varrho }\cap \mathcal{B=\{}%
B_{1},B_{2},B_{5},B_{6}\}.$ This is not a Kurosh-Amitsur radical.
\end{theorem}

\noindent In the next section, we continue on our well-trodden path of
presenting the congruence theory and then the radical theory, in this case
for graphs which do not admit loops. But we will be taken to new
destinations.

\section{Graphs with no loops}

\noindent A congruence theory for graphs that do not allow loops was
introduced by Broere, Heidema and Pretorius [3] and preceeds the congruence
theories for graphs that admit loops and topological spaces in [4] and [30]
respectively. They showed that, as for universal algebra, congruences give
rise to isomorphism theorems as well as the characterization of subdirect
products in terms of the intersection of congruences. When loops are not
allowed, there is a significant restriction on the possible number of
homomorphisms on a graph. A radical theory for such graphs was only
developed recently, and we shall see that this restriction has a
degenerative impact on certain aspects of the radical theory for these
graphs. In particular, the semisimple class of a Hoehnke radical coincides
with the class of all graphs in this case. But in spite of this, there are
non-trivial connectednesses and disconnectednesses for such graphs. It
should be mentioned that an earlier use of a congruence for graphs with no
loops has already appeared in Sabidussi [26], quoting and using results
(unpublished) from his PhD student Fawcett [13]. They used a congruence in a
limited sense as the kernel of a strong graph homomorphism and no general
theory of congruences was mentioned or discussed. In particular, they showed
that subdirect products of graphs can be described in terms of congruences
and they applied these to formulate and prove a version of Birkhoff's
Theorem for graphs: every graph is a subdirect product of subdirectly
irreducible graphs. In their set-up, the subdirectly irreducible graphs are
the complete graphs and the almost complete graphs.

\bigskip

\noindent In this section, all graphs considered are undirected and without
multiple edges. Graphs have non-empty vertex sets, edge sets may be empty
and no loops are allowed. For a graph $G=(V_{G},E_{G}),$ the set of all
possible edges on $G$ will be denoted by $K_{G}:=\{ab\mid a,b\in V_{G},a\neq
b\}.$ Here, for a homomorphism $f;G\rightarrow H$ and $a,b\in V_{G},$ if $%
ab\in E_{G},$ then $f(a)\neq f(b).$

\bigskip

\subsection{Congruences}

\noindent A \textit{congruence on} a graph $G=(V_{G},E_{G})$ [3] is a pair $%
\theta =(\sim ,\mathcal{E})$ such that \newline
(i) $\sim $ is an equivalence relation on $V_{G}$; \newline
(ii) $\mathcal{E}$ is a set of unordered pairs of different elements from $%
V_{G},$ called the \textit{congruence edge set}, with $E_{G}\subseteq 
\mathcal{E\subseteq }$ $K_{G}$; \newline
(iii) when $x\sim y,$ then $xy\not\in \mathcal{E};$ and \newline
(iv) (\textit{Substitution Property} of $\mathcal{E}$ with respect to $\sim
) $ when $x,y,x^{\prime },y^{\prime }\in V_{G}$, $x\sim x^{\prime }$, $y\sim
y^{\prime }$, and $xy\in \mathcal{E}$, then $x^{\prime }y^{\prime }\in 
\mathcal{E}$.\newline
A \textit{strong congruence on }$G$ is a pair $\theta =(\sim ,\mathcal{E})$
where $\sim $ is an equivalence relation on $V_{G},$ $\mathcal{E}=\{xy\mid
x,y\in V_{G}$ and there are $x^{\prime },y^{\prime }\in V_{G}$ with $x\sim
x^{\prime },y\sim y^{\prime }$ and $x^{\prime }y^{\prime }\in E_{G}\}$ and
condition (iii) is fulfilled.

\noindent Requirement (iii) ensures that the equivalence classes $%
[x]:=\{y\in V\mid x\sim y\}$ are independent sets of vertices with respect
to $\mathcal{E}$, i.e., if $a,b\in \lbrack x],$ then $ab\notin E_{G}$. It
can easily be verified that a strong congruence is also a congruence.
Congruences are partially ordered by the relation "contained in": for two
congruences $\alpha =(\sim _{\alpha },\mathcal{E}_{\alpha })$ and $\beta
=(\sim _{\beta },\mathcal{E}_{\beta })$ on $G,$ $\alpha $ is \textit{%
contained in }$\beta ,$ written as $\alpha \subseteq \beta ,$ if $\sim
_{\alpha }\subseteq \sim _{\beta }$ and $\mathcal{E}_{\alpha }\subseteq 
\mathcal{E}_{\beta }.$ We will always use $\Bumpeq $ to denote the identity
relation (diagonal) on a set. The congruence $\iota _{G}:=(\Bumpeq ,E_{G})$
on $G,$ called the \textit{identity congruence} on $G,$ is the smallest
congruence on $G.$ It is a strong congruence on $G$. If $G$ is a complete
graph, i.e., $E_{G}=K_{G}$, then $G$ can have only one congruence namely the
identity congruence $\iota _{G}$.

\bigskip

\noindent \textbf{The kernel of a homomorphism. }Given any graph
homomorphism $f:G\longrightarrow H$, a congruence on $G$, called the \textit{%
kernel of} $f$ and written as $\ker f=(\sim _{f},\mathcal{E}_{f}),$ is
defined by $\sim _{f}=\{(x,y)\mid x,y\in V_{G},f(x)=f(y)\}$ and $\mathcal{E}%
_{f}=\{uv\mid u,v\in V_{G},f(u)f(v)\in E_{H}\}.$ It is immediately clear
that $\ker f$ is a congruence on $G$. With $f$ is also associated the 
\textit{strong kernel of }$f,$ written as \textit{sker} $f=(\sim _{f},%
\mathcal{E}_{sf})$ with the same equivalence relation but $\mathcal{E}%
_{sf}=\{xy\mid x,y\in V_{G}$ and there are $x^{\prime },y^{\prime }\in V_{G}$
with $x\sim _{f}x^{\prime }$, $y\sim _{f}y^{\prime }$ and $x^{\prime
}y^{\prime }\in E_{G}\}.$ This is a strong congruence on $G$ and \textit{sker%
} $f\subseteq \ker f;$ in fact, if $\theta =(\sim _{f},\mathcal{E})$ is any
congruence on $G$ for some $\mathcal{E}$, then \textit{sker} $f\subseteq
\theta .$ If $f$ is a strong homomorphism, then $\ker f=$ \textit{sker} $f.$
Note that a homomorphism $f$ is injective if and only if $\sim _{f}=$ $%
\Bumpeq .$ Moreover, if $f$ is a surjective strong homomorphism, then $f$ is
an isomorphism if and only if $\ker f=\iota _{G}.$

\bigskip

\noindent \textbf{Quotients. }Given any congruence $\theta =(\sim ,\mathcal{E%
})$ on a graph $G=(V_{G},E_{G})$, a new graph, denoted by $G/\theta
=(V_{G/\theta },E_{G/\theta })$ and called \textit{the quotient of }$G$%
\textit{\ by }$\theta ,$ is defined by taking $V_{G/\theta }:=\{[x]\mid x\in
V_{G}\}$ and $E_{G/\theta }:=\{[x][y]\mid xy\in \mathcal{E}\}.$ The natural
(canonical) mapping $p_{\theta }:G\rightarrow G/\theta $ given by $p_{\theta
}(x)=[x]$ is a surjective homomorphism with $\ker p_{\theta }=\theta .$ In
particular, for $\theta =\iota _{G}$ we have $G/\iota _{G}$ is isomorphic to 
$G.$ If $\theta $ is a strong congruence, then $p_{\theta }$ is a strong
homomorphism with \textit{sker} $p_{\theta }=\theta =\ker \theta .$ In
general, if we take $\theta =(\Bumpeq ,\mathcal{E})$ for some suitable $%
\mathcal{E}$ to make $\theta $ a congruence on $G,$ then $G/\theta $ is the
graph with vertex set $V_{G/\theta }=V_{G}$ and edge set $E_{G/\theta }=%
\mathcal{E}$ (here we identify $[x]=\{x\}$ with $x).$

\bigskip

\noindent \textbf{The semilattice of congruences. }For a given graph $G,$ we
denote the set of all congruences on $G$ by $Con(G).$ We already know that $%
Con(G)$ is a partially ordered set with respect to containment $\subseteq .$
But we can say more. Any collection of congruences $\{\theta _{i}=(\sim _{i},%
\mathcal{E}_{i})\mid i\in I\}\subseteq Con(G)$ has a greatest lower bound in 
$Con(G)$ given by $\bigcap\limits_{i\in I}\theta _{i}=(\sim ,\mathcal{E})$
where $a\sim b\Leftrightarrow a\sim _{i}b$ for all $i\in I$ and $ab\in 
\mathcal{E}\Leftrightarrow ab\in \mathcal{E}_{i}$ for all $i\in I.$ This
ensures that $Con(G)$ is a complete meet-semilattice with the meet given by
the intersection as defined above.

\bigskip

\noindent \textbf{Image of a congruence. }Let $f:G\rightarrow H$ be a
homomorphism and $\theta =(\sim ,\mathcal{E})$ a congruence on $G.$ Then $%
f(\theta )$ means the pair $(f(\sim ),f(\mathcal{E}))$ with $f(\sim
):=\{(f(a),f(b))\mid a,b\in V_{G},a\sim b\}\subseteq V_{H}\times V_{H}$ and $%
f(\mathcal{E}):=\{f(a)f(b)\mid ab\in \mathcal{E\}}\subseteq \{xy\mid x,y\in
V_{H}\}.$ Note that $f(\theta )$ need not be a congruence on the graph $H.$
Despite this, for a congruence $\beta =(\sim _{\beta },\mathcal{E}_{\beta })$
on $H,$ we will compare $f(\theta )$ with $\beta $ in the usual sense: $%
f(\theta )\subseteq \beta $ if $f(\sim )\subseteq \sim _{\beta }$ and $f(%
\mathcal{E})\subseteq \mathcal{E}_{\beta }$.

\bigskip

\noindent \textbf{Isomorphism theorems for congruences. }When dealing with
radicals, the basic tools are the appropriate versions of the algebraic
isomorphism theorems for graph congruences. These theorems are discussed in
detail in Broere, Heidema and Pretorius [3], but summarized below for ease
of reference.

\noindent (1) (\textit{First Isomorphism Theorem}) [3] Let $f:G\rightarrow H$
be a surjective homomorphism. Then $G/\ker f$ is isomorphic to $H$.

\bigskip

\noindent The Second Isomorphism has also been given in [3], but we will
present it here using a different (and more suggestive) notation. Let $G$ be
a graph with induced subgraph $H.$ Then a congruence $\theta =(\sim ,%
\mathcal{E})$ on $G$ induces a congruence $H\cap \theta =(\sim _{H},\mathcal{%
E}_{H})$ on $H$ with $\sim _{H}=\{(a,b)\mid a,b\in V_{H}$ and $a\sim b\}$
and $\mathcal{E}_{H}=\{ab\mid a,b\in V_{H}$ with $ab\in \mathcal{E}\}.$ The
mapping $f:H\rightarrow G/\theta $ defined by $f(a)=[a]$ for all $a\in H$ is
a homomorphism with $\ker f=H\cap \theta .$ Now $f(V_{H})$ is a set of
vertices of $G/\theta $ on which we form the induced subgraph of $G/\theta $%
, denoted by $(H+\theta )/\theta .$ Then, by the First Isomorphism Theorem,
we have:

\noindent (2) (\textit{Second Isomorphism Theorem}) [3] Let $H$ be an
induced subgraph of a graph $G.$ Let $\theta $ be a congruence on $G.$ Then $%
H/H\cap \theta \cong (H+\theta )/\theta $ where the latter graph is the
induced subgraph of $G/\theta $ on the vertex set $\{[a]\mid a\in V_{H}\}.$

\bigskip

\noindent (3) (\textit{Third Isomorphism Theorem}) [3] Let $G$ be a graph
with $\theta _{1}=(\sim _{1},\mathcal{E}_{1})$ and $\theta _{2}=(\sim _{2},%
\mathcal{E}_{2})$ two congruences on $G$ for which $\theta _{1}\subseteq
\theta _{2}.$ Then $\theta _{2}/\theta _{1}:=(\sim ,\mathcal{E})$ is a
congruence on $G/\theta _{1}$ where $[a]_{1}\sim \lbrack
b]_{1}\Leftrightarrow a\sim _{2}b$ and $[a]_{1}[b]_{1}\in \mathcal{E}%
\Leftrightarrow ab\in \mathcal{E}_{2}.$ Moreover, $(G/\theta _{1})/(\theta
_{2}/\theta _{1})$ is isomorphic to $G/\theta _{2}.$

\bigskip

\noindent A related result often used, also from [3], is:

\noindent (4) Let $G$ be a graph with $\theta $ a fixed congruence on $G.$
Any congruence $\xi $ of the graph $G/\theta $ is of the form $\alpha
/\theta $ for some congruence $\alpha $ on $G$ with $\theta \subseteq \alpha
.$ Moreover, there is a one-to-one correspondence between $\{\alpha \mid
\alpha $ is a congruence on $G$ with $\theta \subseteq \alpha \}$ and $%
Con(G/\theta )$ which preserves inclusions and intersections.

\bigskip

\noindent \textbf{Products, subdirect products and Birkhoff's Theorem. }For
an index set $I,$ let $G_{i}=(V_{i},E_{i})$ be a graph for all $i\in I.$ The
product $\prod\limits_{i\in I}G_{i}$ of the graphs $G_{i}$ is the graph $%
\prod\limits_{i\in I}G_{i}=(\prod\limits_{i\in I}V_{i},E)$ where $%
\prod\limits_{i\in I}V_{i}$ is just the usual Cartesian product of the sets $%
V_{i}$ and $E=\{fg\mid f,g\in \prod\limits_{i\in I}V_{i}$ with $f(i)g(i)\in
E_{i}$ for all $i\in I\}.$ For every $j\in I,$ the $j$-th \textit{projection}
$\pi _{j}:\prod\limits_{i\in I}G_{i}\rightarrow G_{j}$ defined by $\pi
_{j}(f)=f(j)$ for all $f\in \prod\limits_{i\in I}V_{i}$ is a surjective
homomorphism. An induced subgraph $H$ of $\prod\limits_{i\in I}G_{i}$ is
called a \textit{subdirect product of the graphs} $G_{i},i\in I,$ provided
the restriction of each projection $\pi _{j}$ to $H$ is surjective. As in
universal algebra, subdirect products can be expressed in terms of
congruences and quotients:

\begin{theorem}
\lbrack 3] A graph $G$ is a subdirect product of graphs $G_{i},i\in I,$ if
and only if for every $i\in I$ there are congruences $\theta _{i}$ on $G$
with $G_{i}$ isomorphic to $G/\theta _{i}$ and $\bigcap\limits_{i\in
I}\theta _{i}=\iota _{G}.$
\end{theorem}

\bigskip

\noindent A graph $G$ which has the property that whenever it is a subdirect
product of graphs $G_{i},i\in I,$ then at least one of the $G_{i}$'s must be
isomorphic to $G$ is called \textit{subdirectly irreducible}. In view of the
theorem above, a graph $G$ is subdirectly irreducible if and only if any
intersection of congruences on $G$ which is the identity congruence, must
already include one congruence which is the identity. A graph which is not
subdirectly irreducible, is called \textit{subdirectly reducible}. The
corresponding notions for algebra are important, especially in the context
of the well-known theorem of Birkhoff [2]: every non-trivial algebra is a
subdirect product of subdirectly irreducible algebras. And, of course,
knowing exactly what the subdirectly irreducible algebras are is the essence
of this result. This theorem has meaning and validity for many other
mathematical structures as well. For example, every non-trivial topological
spaces is a subdirect product of copies of the Sierspi\'{n}ski space and the
two-element indiscrete space (see, for example, Proposition 2.4 in [31]),
noting that both these two two-element spaces are subdirectly irreducible
topological spaces. Birkhoff's Theorem is also valid in the category of
graphs that admit loops [35], where it says that any non-trivial graph is a
subdirect product of subdirectly irreducible graphs. Here the subdirectly
irreducible graphs are $B_{4},$ $B_{5},$ $B_{6}$ (see Section 4.3 for the
definitions of the $B_{i}$'s) and $A_{3}$ where $A_{3}$ is the three-vertex
graph $V_{A_{3}}=\{0,1,2\}$ and $E_{A_{3}}=\{00,11,22,01,21\}.$ Sabidussi
and Fawcett's version of Birkhoff's theorem (cf. [26] and [13]) is not valid
in our more general setting. In our case, a graph is subdirectly irreducible
if and only if it is a complete graph and then:

\begin{theorem}
\lbrack 36] Every graph is a subdirect product of complete graphs.
\end{theorem}

\subsection{Radical theory}

\noindent The prohibition on loops imposes a significant restriction on the
admissible maps between graphs and will consequently have an impact on the
radical theory for these graphs. In fact, it will be shown that the radical
theory in the category of graphs that do not admit loops has significant
differences with all existing radical theories. In particular, all Hoehnke
radicals are degenerative, but there are non-trivial connectednesses
(radical classes) and disconnectednesses (semisimple classes). Contrary to
radical theories in other categories, this means that objects need not have
maximal semisimple images. We start by defining a Hoehnke radical and report
on their salient features. But this is really just a wild goose chase, since
they are all trivial. Then we define the connectednesses and
disconnectednesses, give non-trivial examples and show that they always come
as complementary pairs. More detail on what is presented here can be found
in [36].

\bigskip

\noindent We will work in a universal class $\mathcal{W}$ of graphs
(non-empty, closed under homomorphic images and closed under the taking of
subgraphs (= strongly hereditary)). From the definition, it follows that $%
\mathcal{W}$ contains a one-vertex graph, and consequently all one-vertex
graphs. We identify all the one-vertex-graphs, denote them by $T$ and call
them the \textit{trivial graphs}. For any graph $G,$ the \textit{completion
of }$G$ is the graph $G^{c}$ with the same vertex set as $G$ and edge set $%
K_{G}.$ Since $G^{c}$ is a homomorphic image of $G,$ $\mathcal{W}$ contains
the completion of all graphs $G\in \mathcal{W.}$ We will assume $\mathcal{W}$
is non-trivial, i.e., it has at least one graph with two or more vertices.
Since $\mathcal{W}$ is strongly hereditary, it thus contains the two-vertex
graph with no edges $B_{1}$ and its completion $B_{1}^{C}=B_{2}.$ All
considerations relating to the radicals of graphs will be inside the class $%
\mathcal{W}.$

\begin{definition}
An $H$-\textit{radical} on $\mathcal{W}$ is a function $\rho $ that assigns
to every graph $G$ in $\mathcal{W}$ a congruence $\rho (G)=\rho _{G}$ on $G$
such that: \newline
(H1) if $f:G\rightarrow H$ is a surjective homomorphism in $\mathcal{W},$
then $f(\rho _{G})\subseteq \rho _{H}$; and\newline
(H2) for all graphs $G\in \mathcal{W},$ $\rho (G/\rho _{G})=\iota _{G/\rho
_{G}}$, the identity congruence on $G/\rho _{G}$. \newline
\end{definition}

\noindent For an $H$-radical $\rho ,$ the class $\mathcal{S}_{\rho }=\{G\in 
\mathcal{W}\mid \rho (G)=\iota _{G}\}$ is called the associated \textit{%
semisimple class} and $\mathcal{R}_{\rho }=\{G\in \mathcal{W}\mid G/\rho
_{G} $ is the one-element graph\} is the associated \textit{radical class}. $%
\mathcal{R}_{\rho }$ is always \textit{homomorphically closed}, i.e., if $%
G\in \mathcal{R}_{\rho }$ and $f:G\rightarrow H$ is a surjective
homomorphism, then $H\in \mathcal{R}_{\rho }$ and $\mathcal{S}_{\rho }\cap 
\mathcal{R}_{\rho }=\{T\}.$ The essence of this radical is given in the next
result. A class of graphs $\mathcal{M}$ in $\mathcal{W}$ is an \textit{%
abstract class} provided it contains all the one element graphs in $\mathcal{%
W}$ and it is closed under isomorphic copies. All subclasses of $\mathcal{W}$
under consideration will be assumed to be abstract, even though it may not
always be explicitly stated. For a class $\mathcal{M}$ in $\mathcal{W},$ we
use $\overline{\mathcal{M}}$ to denote the \textit{subdirect closure} of $%
\mathcal{M},$ i.e., the class of all graphs that are subdirect products of
graphs from $\mathcal{M}.$ Clearly $\mathcal{M}\subseteq \overline{\mathcal{M%
}}$ and we say $\mathcal{M}$ is \textit{subdirectly closed} if $\mathcal{M}=%
\overline{\mathcal{M}}.$

\bigskip

\begin{theorem}
(1) Let $\rho $ be a mapping that assigns to any graph $G$ in $\mathcal{W}$
a congruence $\rho (G)=\rho _{G}$ on $G.$ If $\rho $ is an $H$-radical on $%
\mathcal{W,}$ then there is an abstract class of graphs $\mathcal{M}$ in $%
\mathcal{W}$ such that for all $G$ in $\mathcal{W},$ $\rho (G)=\cap \{\theta
\mid \theta $ is a congruence on $G$ for which $G/\theta \in \mathcal{M}\}.$
Furthermore, $\mathcal{S}_{\rho }=$ $\overline{\mathcal{M}}$ and $\mathcal{S}%
_{\rho }$ is closed under subdirect products.

\noindent (2) For any abstract class of graphs $\mathcal{M}$ in $\mathcal{W}$
for which every $G$ in $\mathcal{W}$ has a congruence $\theta $ on $G$ with $%
G/\theta \in \mathcal{M}$, define a mapping $\rho $ by $\rho (G)=\cap
\{\theta \mid \theta \in Con(G)$ with $G/\theta \in \mathcal{M}\}.$ Then $%
\rho $ is an $H$- radical and $\mathcal{S}_{\rho }=$ $\overline{\mathcal{M}}$%
.
\end{theorem}

\noindent In the context of the theorem above, we say that the $H$-radical $%
\rho $ is determined by the class $\mathcal{M}$ if $\mathcal{M}\subseteq 
\mathcal{S}_{\rho }$ and $\overline{\mathcal{M}}=\mathcal{S}_{\rho }$; the
class $\mathcal{M}$ is not necessarily unique. In categories where
congruences can be identified by a distinguished subobject (e.g. for rings,
a congruence is completely determined by an ideal which is the congruence
class of the additive identity), any class $\mathcal{M}$ of objects will
fulfil the requirement of (2) above, but in general it need not be the case.
Fortunately we can say exactly when it will be in our case. Let $\mathcal{K}%
_{\mathcal{W}}$ be the set of all complete graphs in $\mathcal{W}.$

\begin{lemma}
Let $\mathcal{M}$ be an abstract class of graphs in $\mathcal{W}.$ Then for
every $G$ in $\mathcal{W}$ there is a congruence $\theta $ on $G$ with $%
G/\theta \in \mathcal{M}$ if and only if $\mathcal{K}_{\mathcal{W}}\subseteq 
\mathcal{M}.$
\end{lemma}

\noindent This result is not good for the well-being of Hoehnke radicals in
this universal class, for we have:

\begin{theorem}
Let $\rho $ be an $H$-radical on $\mathcal{W.}$ Then $\mathcal{S}_{\rho }=%
\mathcal{W},$ i.e., $\rho (G)=\iota _{G}$ for all $G\in \mathcal{W},$ and $%
\mathcal{R}_{\rho }=\{T\}.$
\end{theorem}

\noindent Any $H$-radical on a universal class of graphs that do not admit
loops is thus degenerative in the sense that $\rho (G)=\iota _{G}$ for all $%
G.$ On the other hand, in general radical theory, the semisimple objects are
usually regarded as well-behaved and sought after objects. One
interpretation of the result $\mathcal{S}_{\rho }=\mathcal{W}$ is that all
the graphs that do not admit loops are good graphs! Next we define the
connectednesses and disconnectednesses in the universal class $\mathcal{W}.$
A class $\mathcal{C}\subseteq \mathcal{W}$ is a \textit{connectedness }(%
\textit{KA-radical class}) if it satisfies the following condition: A graph $%
G$ is in $\mathcal{C}$ if and only if every non-trivial homomorphic image of 
$G$ has a non-trivial induced subgraph which is in $\mathcal{C}.$ A class $%
\mathcal{D}\subseteq \mathcal{W}$ is a \textit{disconnectedness
(KA-semisimple class)} if it satisfies the following condition: A graph $G$
is in $\mathcal{D}$ if and only if every non-trivial induced subgraph of $G$
has a non-trivial homomorphic image which is in $\mathcal{D}.$ The trivial
graph $T$ is always in any connectedness and also in any disconnectedness.
If $\mathcal{M}\subseteq \mathcal{W}$ is a hereditary class, then $\mathcal{%
UM}:=\{G\in \mathcal{W\mid G}$ has no non-trivial homomorphic image in $%
\mathcal{M}\}$ is a connectedness and if $\mathcal{H}\subseteq \mathcal{W}$
is a homomorphically closed class, then $\mathcal{SH}:=\{G\in \mathcal{W}%
\mid G$ has no non-trivial induced subgraph in $\mathcal{H}\}$ is a
disconnectedness.

\begin{proposition}
Any connectedness is homomorphically closed and any disconnectedness is
strongly hereditary and closed under subdirect products.
\end{proposition}

\noindent From the preceding, we thus have: If $\mathcal{C}$ is a
connectedness, then $\mathcal{D}:=\mathcal{SC}$ is a disconnectedness,
called the \textit{disconnectedness corresponding to} $\mathcal{C},$ and if $%
\mathcal{D}$ is a disconnectedness, then $\mathcal{C}:=\mathcal{UD}$ is a
connectedness, called the \textit{connectedness corresponding to} $\mathcal{D%
}.$ Moreover, it can be shown that a class $\mathcal{C}\subseteq \mathcal{W}$
is a connectedness if and only if $\mathcal{C}=\mathcal{USC}$ and a class $%
\mathcal{D}\subseteq \mathcal{W}$ is a disconnectedness if and only if $%
\mathcal{D}=\mathcal{SUD}.$ By Theorem 5.6 we immediately have:

\begin{proposition}
Let $\mathcal{D}$ be a disconnectedness. Suppose there is an $H$-radical $%
\rho $ such that its semisimple class $\mathcal{S}_{\rho }=\mathcal{D}.$
Then $\mathcal{D}=\mathcal{W}$ and the corresponding connectedness $\mathcal{%
C}$ is $\mathcal{C}=\{T\}.$
\end{proposition}

\noindent The significance of this result is in stark contrast to the
following feature of all known radical theories, which, for example we have
seen in all three the previous sections. For a radical theory in a given
universal class $\mathcal{W}$ of objects and a semisimple class $\mathcal{D}$
in $\mathcal{W},$ any object in $\mathcal{W}$ has a maximal semisimple image
in the following sense. For any $A\in \mathcal{W},$ there is a surjective
morphism $q:A\rightarrow A_{\mathcal{D}}$ with $A_{\mathcal{D}}\in \mathcal{D%
}$ such that for any surjective morphism $f:A\rightarrow D$ with $D\in 
\mathcal{D},$ there is a morphism $g:A_{\mathcal{D}}\rightarrow D$ for which 
$g\circ q=f.$ If there is a congruence theory available in $\mathcal{W},$
then this means for every $A\in \mathcal{W},$ there is a congruence $\delta
_{A}$ on $A$ with $A/\delta _{A}\cong A_{\mathcal{D}}$ and for every
congruence $\theta $ on $A$ for which $A/\theta \in \mathcal{D},$ we have $%
\theta \subseteq \delta _{A}.$ Thus we have an $H$-radical $\delta $ on $%
\mathcal{W}$ with $\mathcal{S}_{\delta }=\mathcal{D}$ and $\mathcal{R}%
_{\delta }=\mathcal{UD}.$ In view of Proposition 5.8 above, this raises the
question whether there are any non-trivial connectednesses and
disconnectednesses in $\mathcal{W}.$ By example we will show that there are,
and we will give a very special example (which turns out to be not so
special after all). For a connectedness $\mathcal{C}$ with corresponding
disconnectedness $\mathcal{D}$ in $\mathcal{W,}$ we know $\mathcal{C}\cap 
\mathcal{D}=\{T\}$ and the pair $(\mathcal{C},\mathcal{D})$ is called 
\textit{complementary }if $\mathcal{C}\cup \mathcal{D}=\mathcal{W}.$ Of
course, there are two trivial such pairs, namely $(\mathcal{W},\{T\})$ and $%
(\{T\},\mathcal{W}).$ The existence of non-trivial complementary pairs was
first observed in the category of graphs which allow loops by Fried and
Wiegandt [14]. Such pairs are also to be found in the radical theory for $%
\mathcal{S}$-acts, see Wiegandt [38], with a general condition for their
existence given in [37] (which, however, is not applicable here).

\begin{example}
Suppose $\mathcal{K,}$ the class of all complete graphs, is contained in $%
\mathcal{W}.$ The class $\mathcal{K}$ is homomorphically closed, hence $%
\mathcal{D}:=\mathcal{SK}$ is a disconnectedness and the corresponding
connectedness is $\mathcal{C}=\mathcal{USK}.$ It follows that $\mathcal{D}%
=\{G\in \mathcal{W}\mid E_{G}=\emptyset \}$ and $\mathcal{C}=\{G\in \mathcal{%
W}\mid G=T$ or $E_{G}\neq \emptyset \}$ which is clearly a complementary
pair. It is a non-trivial pair since $\mathcal{W}$ contains all the complete
graphs and $\mathcal{W}$ is strongly hereditary (so $\mathcal{D}\neq \{T\}).$
Moreover, no graph $G$ in $\mathcal{W}$ with $E_{G}\neq \emptyset $ has a
homomorphic image in $\mathcal{D}$ and thus certainly not a maximal one. $%
\mathcal{D}$ is not a connectedness (not homomorphically closed) and $%
\mathcal{C}$ is not a disconnectedness (not strongly hereditary). Any graph $%
G\in \mathcal{W}$ has a maximal homomorphic image in $\mathcal{C}$ in the
following sense: the completion $G^{c}$ of $G$ is a homomorphic image of $G,$
$G^{c}$ is in $\mathcal{C}$ and if $H$ is any other complete homomorphic
image of $G,$ then $H\subseteq G^{c}.$
\end{example}

\noindent This example shows a unique feature of the general radical theory
in this universal class. Here is a very natural universal class with a
radical theory in which all the Hoehnke radicals degenerate and it contains
non-trivial KA-radicals. But there are more surprises. Recall, $B_{1}$ is
the graph with two vertices and no edges. For any connectedness $\mathcal{C}$
with corresponding disconnectedness $\mathcal{D},$ we must have $B_{1}\in 
\mathcal{C}$ or $B_{1}\in \mathcal{D}$. By $K_{n}$ we denote the complete
graph on $n$ vertices, $n\geq 1.$ Then $K_{2}=B_{1}^{c}$ and $K_{1}=T.$ The
next result shows that there is really nothing special about the example
above in this universal class.

\bigskip

\begin{proposition}
Let $\mathcal{C}\subseteq \mathcal{W}$ be a connectedness with corresponding
disconnectedness $\mathcal{D}=\mathcal{SC}.$ Then $(\mathcal{C},\mathcal{D})$
is a complementary pair. If $B_{1}\in \mathcal{C},$ then the pair $(\mathcal{%
C},\mathcal{D})=(\mathcal{W},\{T\})$ is trivial, and if $B_{1}\in \mathcal{D}%
,$ then the pair $(\mathcal{C},\mathcal{D})$ need not be trivial.
\end{proposition}

\noindent This is another special feature of the radical theory in this
universal class $\mathcal{W}$: any connectedness (respt. disconnectedness)
comes as a connectedness (respt. disconnectedness) in a complementary pair.
We conclude with a general example, which has Example 5.9 as a special case.

\bigskip

\begin{example}
Suppose $K_{n}\in \mathcal{W}$ for $n=1,2,3,...$ . For each $n\geq 1,$ let $%
\mathcal{C}_{n}:=\{G\in \mathcal{W}\mid G=T$ or if $G\neq T,$ then $K_{n}$
is an induced subgraph of $G\}$ and let $\mathcal{D}_{n}:=\{G\in \mathcal{W}%
\mid G=T$ or if $G\neq T,$ then $K_{n}$ is not an induced subgraph of $G\}.$
Then $(\mathcal{C}_{n},$ $\mathcal{D}_{n})$ forms a complementary pair of
connectednesses and disconnectednesses. Clearly $\mathcal{C}_{1}=\mathcal{W}$
and $\mathcal{C}_{2}=\{G\in \mathcal{W}\mid G=T$ or $E_{G}\neq \emptyset \}.$
Let $\mathcal{C}$ be any connectednesss in $\mathcal{W}.$ If $\mathcal{C\neq
W}$, then $\mathcal{C}\subseteq \mathcal{C}_{2}$. This means $\mathcal{W}$
has a largest proper connectedness $\mathcal{C}_{2}$ and thus a smallest
non-trivial disconnectedness $\mathcal{D}_{2}=\{G\in \mathcal{W}\mid
E_{G}=\emptyset \}.$ If $\mathcal{C}$ contains a graph with $n$ vertices, $%
n\geq 2,$ then $K_{n}\in \mathcal{C}$ in which case $\mathcal{C}%
_{n}\subseteq \mathcal{C}.$ If the universal class $\mathcal{W}$ contains at
least one graph with an infinite number of vertices, then $\mathcal{C}%
:=\{G\in \mathcal{W}\mid $ $G=T$ or $E_{G}$ is infinite\} and $\mathcal{D}%
:=\{G\in \mathcal{W}\mid E_{G}$ is finite\} is a corresponding pair of
non-trivial connectedness and disconnectedness respectively and $\mathcal{C}%
_{n}\nsubseteqq \mathcal{C}$ for all $n\geq 1.$
\end{example}

\bigskip

\noindent \textbf{Acknowledgement. }Thanks are due to Izak Broere who drew
my attention to the existence of a theory of congruences for graphs and also
to Johannes Heidema's idea to define a Hoehnke radical for graphs.

\bigskip

\noindent \textbf{References}

\noindent \lbrack 1] A.V. Arhangel'ski\u{\i} and R. Wiegandt.
Connectednesses and disconectednesses in topology. \textit{General Topology
and Appl.} \textbf{5 }(1975), 9-33.\newline
[2] {G. Birkhoff. Subdirect unions in universal algebra. }\textit{Bull.
Amer. Math. Soc}{.}, \textbf{50 {(}}1944\textbf{), }{764 -- 768.}\newline
[3] I. {Broere, J. Heidema and L.M. Pretorius. }Graph congruences and what
they connote. \textit{Quaestiones Mathematicae} \textbf{41}(8), 1045 - 1059
(2018).\newline
[4] I. Broere, J. Heidema and S. Veldsman. Congruences and Hoehnke Radicals
on Graphs. \textit{Discussiones Mathematicae Graph Theory} \textbf{40}(4)
(2020), 1067-1084.\newline
[5] A. Buys, N.J. Groenewald and S. Veldsman. Radical and semisimple classes
in categories. \textit{Quaestiones Mathematicae} \textbf{4} (1981), 205-220.%
\newline
[6] G. Castellini. Connectedness classes. \textit{Quaestiones Mathematicae} 
\textbf{23} (3) (2000), 313-334.\newline
[7] G. Castellini. Disconnectedness classes. \textit{Quaestiones Mathematicae%
} \textbf{24} (1) (2001), 75-92.\newline
[8] G. Castellini. Connectedness with respect to a closure operator. \textit{%
Appl. Categ. Structures} \textbf{9} (2001), 285-301.\newline
[9] G. Castellini and D. Holgate. A link between two connectedness notions. 
\textit{Appl. Categ. Structures} \textbf{11} (2003), 473-486.\newline
[10] M.M Clementino. On connectedness via closure operators. \textit{Appl.
Categ. Structures} \textbf{9}(2001), 539-556.\newline
[11] M.M. Clementino and W. Tholen. Separation versus connectedness. \textit{%
Topology Appl}. \textbf{75} (1997), 143-181.\newline
[12] S.E. Dickson. A torsion theory for Abelian categories. \textit{Trans.
Amer. Math. Soc}. \textbf{121} (1966), 223-235.\newline
[13] B. Fawcett. \textit{Graphs and ultrapowers}. PhD thesis, McMaster
University, 1969.\newline
[14] E. Fried and R. Wiegandt. Connectednesses and disconnectednesses of
graphs. \textit{Algebra Universalis }\textbf{5 }(1975)\textbf{, }411 -- 428.%
\newline
[15] E. Fried and R. Wiegandt. Abstract relational structures, I (General
Theory). \textit{Algebra Universalis} \textbf{15} (1982), 1-21.\newline
[16] E. Fried and R. Wiegandt. Abstract relational structures, II (Torsion
theory). \textit{Algebra Universalis} \textbf{15} (1982), 22-39.\newline
[17] B.J. Gardner and R. Wiegandt. \textit{Radical Theory of Rings. }Marcell
Dekker Inc., USA,\textit{\ }2004.\newline
[18] {H.-J. Hoehnke. }Radikale in allgemeinen Algebren. \textit{%
Mathematische Nachrichten} \textbf{32} (1966) 347--383.\newline
[19] M. Holcombe and R. Walker. Radicals in categories. \textit{Proc.
Edinburgh Math. Soc}. \textbf{21} (1978/9), 111-128.\newline
[20] L. M\'{a}rki, R. Mlitz and R. Wiegandt. A general Kurosh-Amitsur
radical theory. \textit{Comm. Algebra} \textbf{16} (1988), 249--305.\newline
[21] R. Mlitz. Radical and semisimple classes of $\Omega $-groups. \textit{%
Proc. Edinburgh Math. Soc.} \textbf{23} (1980), 37-41.\newline
[22] R. Mlitz and S. Veldsman. Radicals and subdirect decompositions of $%
\Omega $-groups. \textit{J. Austal. Math. So}c. \textbf{48} (1990), 171-198.%
\newline
[23] G. Preu\ss . Trennung und Zusammenhang. \textit{Monatsh. Math. }\textbf{%
74} (1970), 70-87.\newline
[24] G. Preu\ss . E-zusammenh\"{a}ngende R\"{a}ume. \textit{Manuscripta Math}%
. \textbf{3} (1970), 331-342.\newline
[25] G. Preu\ss . Eine Galois-Korrespondenz in der Topologie. \textit{%
Monatsh. Math}. \textbf{75 }(1971), 447-452.\newline
[26] G. Sabidussi. Subdirect representations of graphs. \textit{Coll. Math.
Soc. Janos Bolyai 10 Infinite and finite sets, Kesztheley (Hungary, 1973)}.
Vol III. A. Hajnal e.a. (reds), 1199 - 1226. North Holland, 1974. \newline
[27] E.G. \v{S}ul'ge\u{\i}fer. On the general theory of radicals in
categories. \textit{Amer. Math. Soc. Transl}. \textbf{59} (1966), 150-162.%
\newline
[28] A. Suli\'{n}ski. The Brown-McCoy radical in categories. \textit{Fund.
Math}. \textbf{59} (1966), 23-41.\newline
[29] S. Veldsman. On the characterization of radical and semisimple classes
in categories. \textit{Comm. Algebra} \textbf{10} (1982), 913-938.\newline
[30] S. Veldsman. Congruences on topological spaces with an application to
radical theory. \textit{Algebra Universalis} \textbf{80}(2019), article 25.%
\newline
[31] S. Veldsman. Hereditary torsion theories and connectednesses and
disconnectednesses for topological spaces. \textit{East-West J. Math.} 
\textbf{21}(2019), 85-102.\newline
[32] S. Veldsman. Topological connectednesses and congruences. \textit{%
Quaestiones Mathematicae} \textbf{44}(12) 2021, 1757-1772.\newline
[33] S. Veldsman. Hereditary torsion theories for graphs. \textit{Acta
Mathematica Hungarica} \textbf{163}(2), (2021), 363-378.\newline
[34] S. Veldsman. Connectednesses of graphs and congruences. \textit{%
Asian-European Journal of Mathematics} \textbf{14} (10) (2021), 19 pages.%
\newline
[35] S. Veldsman. Congruences and subdirect representations of graphs.%
\textit{\ Electronic Journal of Graph Theory and Applications} \textbf{8}
(2020), 123-132.\newline
[36] S. Veldsman. A radical theory for graphs that do not admit loops. 
\textit{Periodica Mathematica Hungarica}.
https://doi.org/10.1007/s10998-022-00496-0. Published online Aug 2022.%
\newline
[37] S. Veldsman and R. Wiegandt. On the existence and non-existence of
complementary radical and semisimple classes. \textit{Quaestiones
Mathematicae} \textbf{7} (1984), 213 -- 224.\newline
[38] R. Wiegandt. Radical and Torsion Theory for Acts. \textit{Semigroup
Forum} \textbf{72} (2006), 312-328.

\bigskip

\noindent \textbf{Stefan Veldsman}

Nelson Mandela University, Port Elizabeth, SOUTH AFRICA

and

La Trobe University, Melbourne, AUSTRALIA

Email: veldsman@outlook.com

\end{document}